\documentclass{gtart_a}
\pdfoutput=1
\usepackage{pinlabel}

\title[Gromov--Witten invariants of configurations of curves]{The local Gromov--Witten invariants of\\configurations of rational curves}

\author[Karp]{Dagan Karp}
\givenname{Dagan}
\surname{Karp}
\address{Department of Mathematics\\University of California at
Berkeley\\\newline
California 94720-3840\\USA}
\email{dkarp@math.berkeley.edu}
\urladdr{}

\author[Liu]{Chiu-Chu Melissa Liu}
\givenname{Chiu-Chu Melissa}
\surname{Liu}
\address{Department of Mathematics\\Northwestern University\\\newline
Evanston\\Illinois 60208-2370\\USA}
\email{ccliu@math.northwestern.edu}
\urladdr{}

\author[Mari\~no]{Marcos Mari\~no}
\givenname{Marcos}
\surname{Mari\~no}
\address{Department of Physics\\CERN\\\newline
Geneva 23\\CH-1211\\Switzerland}
\email{marcos@mail.cern.ch}
\urladdr{}

\volumenumber{10}
\issuenumber{}
\publicationyear{2006}
\papernumber{5}
\lognumber{0646}
\startpage{115}
\endpage{168}

\doi{}
\MR{}
\Zbl{}


\received{23 August 2005}
\revised{23 January 2006}
\accepted{25 November 2005}
\publishedonline{7 March 2006}
\published{7 March 2006}
\proposed{Jim Bryan}
\seconded{Paul Goerss, Eleny Ionel}
\corresponding{}
\editor{}
\version{}

\begin{htmlabstract}
We compute the local Gromov&ndash;Witten invariants of certain configurations
of rational curves in a Calabi&ndash;Yau threefold. These configurations
are connected subcurves of the "minimal trivalent configuration",
which is a particular tree of <b>P</b><sup>1</sup>'s with specified formal
neighborhood. We show that these local invariants are equal to certain
global or ordinary Gromov&ndash;Witten invariants of a blowup of <b>P</b><sup>3</sup>
at points, and we compute these ordinary invariants using the geometry of
the Cremona transform. We also realize the configurations in question as
formal toric schemes and compute their formal Gromov&ndash;Witten invariants
using the mathematical and physical theories of the topological vertex. In
particular, we provide further evidence equating the vertex amplitudes
derived from physical and mathematical theories of the topological vertex.
\end{htmlabstract}

\begin{asciiabstract}
We compute the local Gromov-Witten invariants of certain
configurations of rational curves in a Calabi-Yau threefold. These
configurations are connected subcurves of the `minimal trivalent
configuration', which is a particular tree of P^1's with specified
formal neighborhood. We show that these local invariants are equal to
certain global or ordinary Gromov-Witten invariants of a blowup of
P^3 at points, and we compute these ordinary invariants using the
geometry of the Cremona transform. We also realize the configurations
in question as formal toric schemes and compute their formal
Gromov--Witten invariants using the mathematical and physical theories
of the topological vertex. In particular, we provide further evidence
equating the vertex amplitudes derived from physical and mathematical
theories of the topological vertex.
\\
We compute the local Gromov--Witten invariants of certain configurations
of rational curves in a Calabi--Yau threefold. These configurations
are connected subcurves of the ``minimal trivalent configuration'',
which is a particular tree of $\mathbb{P}^1$'s with specified formal
neighborhood. We show that these local invariants are equal to certain
global or ordinary Gromov--Witten invariants of a blowup of $\mathbb{P}^3$
at points, and we compute these ordinary invariants using the geometry of
the Cremona transform. We also realize the configurations in question as
formal toric schemes and compute their formal Gromov--Witten invariants
using the mathematical and physical theories of the topological vertex. In
particular, we provide further evidence equating the vertex amplitudes
derived from physical and mathematical theories of the topological vertex.
\end{asciiabstract}

\newcommand{\cM}{\mathcal{M}}

\newcommand{\cW}{\mathcal{W}}

\newcommand{\cO}{\mathcal{O}}

\newcommand{\Aut}{\mathrm{Aut}}
\newcommand{\Hom}{\mathrm{Hom}}
\newcommand{\vir}{{\mathrm{vir} }}

\newcommand{\bC}{\mathbb{C}}
\newcommand{\bP}{\mathbb{P}}  
\newcommand{\QQ}{\mathbb{Q}}
\newcommand{\bQ}{\mathbb{Q}}
\newcommand{\bS}{\mathbb{S}}
\newcommand{\bT}{\mathbb{T}}

\newcommand{\bZ}{\mathbb{Z}}

\newcommand{\bp}{\mathbf{p}}
\newcommand{\bw}{\mathbf{w}}
\newcommand{\vn}{\mathbf{n}}
\newcommand{\vp}{\mathbf{p}}

\newcommand{\vo}{\mathbf{0}}
\newcommand{\vd}{\mathbf{d}}     
\newcommand{\vt}{\mathbf{t}}

\newcommand{\hY}{\hat{Y}}

\newcommand{\lam}{\lambda}
\newcommand{\la}{\lambda}
\newcommand{\si}{\sigma}

\newcommand{\tF}{\tilde{F}}
\newcommand{\tN}{\tilde{N}}
\newcommand{\tX}{\tilde{X}}
\newcommand{\tZ}{\tilde{Z}}
\newcommand{\vmu}{{\vec{\mu}}}
\newcommand{\vnu}{{\vec{\nu}}}

\newcommand{\bu}{\bullet}
\newcommand{\Mbar}{{}\mskip3mu\overline{\mskip-3mu\cM\mskip-1mu}\mskip1mu}
\newcommand{\up}[1]{ {{#1}^1,{#1}^2,{#1}^3} }
\newcommand{\lo}[1]{ {{#1}_1,{#1}_2,{#1}_3} }
\newcommand{\three}[1]{ {(\lo{#1})} }
\newcommand{\braket}{\left< \; \right>}
\newcommand{\im}{\operatorname{Im}}
\newcommand{\be}{\begin{equation}}
\newcommand{\ee}{\end{equation}}

\AtBeginDocument{\let\bar\wbar\let\tilde\wtilde\let\hat\what}
\def\cnewtheorem#1[#2]#3{\newtheorem{#1}{#3}}

\newtheorem{prop}{Proposition}
\newtheorem{fact}{Fact}
\cnewtheorem{theo}[prop]{Theorem}
\cnewtheorem{lemm}[prop]{Lemma}
\cnewtheorem{coro}[prop]{Corollary}
\cnewtheorem{rema}[prop]{Remark}
\cnewtheorem{exam}[prop]{Example}
\cnewtheorem{defi}[prop]{Definition}
\cnewtheorem{conj}[prop]{Conjecture}

\makeatletter
\let\c@theo\c@prop
\let\c@lemm\c@prop
\let\c@coro\c@prop
\let\c@rema\c@prop
\let\c@exam\c@prop
\let\c@defi\c@prop
\let\c@conj\c@prop
\makeatother

\begin{document}

\begin{abstract}
We compute the local Gromov--Witten invariants of certain configurations
of rational curves in a Calabi--Yau threefold. These configurations
are connected subcurves of the ``minimal trivalent configuration'',
which is a particular tree of $\mathbb{P}^1$'s with specified formal
neighborhood. We show that these local invariants are equal to certain
global or ordinary Gromov--Witten invariants of a blowup of $\mathbb{P}^3$
at points, and we compute these ordinary invariants using the geometry of
the Cremona transform. We also realize the configurations in question as
formal toric schemes and compute their formal Gromov--Witten invariants
using the mathematical and physical theories of the topological vertex. In
particular, we provide further evidence equating the vertex amplitudes
derived from physical and mathematical theories of the topological vertex.
\end{abstract}

\maketitle

\section{Introduction}\label{sec:introduction}

Let $Z$ be a closed subvariety of a smooth projective threefold $X$
such that $X$ is a local Calabi--Yau threefold near $Z$.  In some cases,
the contribution to the Gromov--Witten invariants of $X$ by maps to $Z$
can be isolated and defines local Gromov--Witten invariants of $Z$ in
$X$. Information obtained from the study of local Gromov--Witten theory
can be used to gain insight into Gromov--Witten theory in general. This
has led to a great amount of interest in the subject.

The study of the local invariants of curves in a Calabi--Yau threefold
has a particularly rich history.  Their study goes back to the famous
Aspinwall--Morrison formula for the local invariants of a single  $\bP^1$
smoothly embedded in a Calabi--Yau threefold with normal bundle $\cO
(-1) \oplus \cO (-1)$; this result is studied by Aspinwall and Morrison
\cite{AM}, Cox and Katz \cite{CK}, Faber and Pandharipande \cite{FP},
Kontsevich \cite{Ko}, Lian, Liu and Yau \cite{LLY}, Manin \cite{Ma},
Pandharipande \cite{P}, and Voisin \cite{V}. The local invariants of
nonsingular curves of any genus have been completely determined by Bryan
and Pandharipande \cite{BP01,BP03,BP04}. In \cite{BKaL}, Bryan, Katz and
Leung computed local invariants of certain rational curves with nodal
singularities, and in particular, contractible $ADE$ configurations of
rational curves. The local invariants of the closed topological vertex,
which is a configuration of three $\bP^1$'s meeting in a single triple
point, were computed by Bryan and the first author \cite{BK}.

In this paper, we will compute local invariants of certain configurations of 
rational curves. The configurations considered in this paper are all connected
subtrees of the {\em minimal trivalent configuration}, which is a configuration of
three chains of $\bP^1$'s meeting in a triple point (see Figure 1 below).
A precise description of the formal neighborhood will be given in \fullref{sec:construction}.
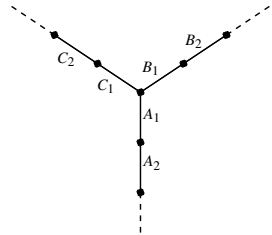
\begin{figure}[ht!]
\labellist\tiny
\hair=1pt
\pinlabel {$A_1$} [l] at 165 153
\pinlabel {$A_2$} [l] at 165 93
\pinlabel {$B_1$} [br] at 191 201
\pinlabel {$B_2$} [br] at 245 235
\pinlabel {$C_1$} [tr] at 135 201
\pinlabel {$C_2$} [tr] at 85 235
\endlabellist
\centerline{\includegraphics[scale=0.3]{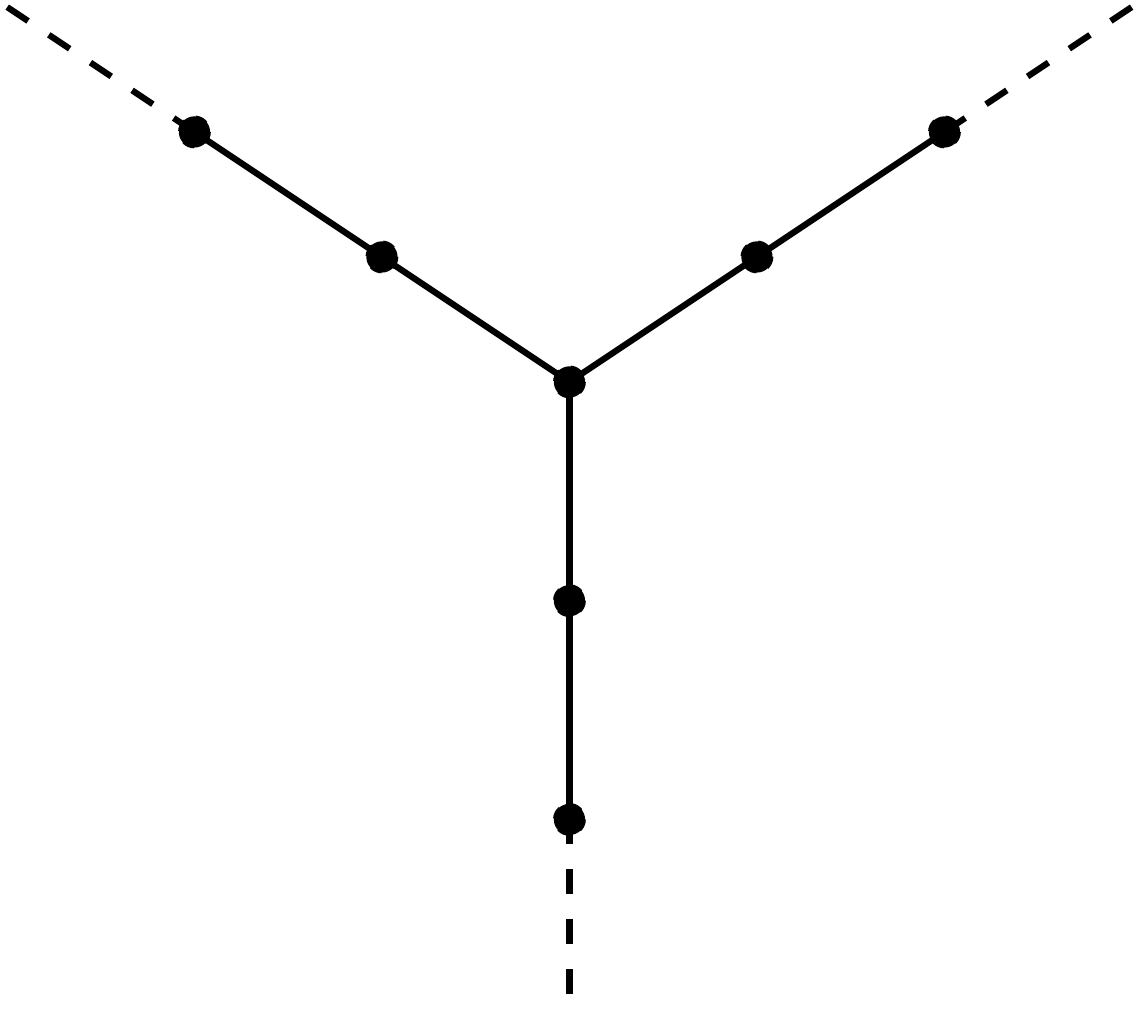}}
\caption{The minimal trivalent configuration $Y^N=\bigcup_{i=1}^N
A_i\cup B_i \cup C_i$.  The normal bundles of $A_1,B_1,C_1$ are isomorphic
to $\cO(-1)\oplus \cO(-1)$; the normal bundle of any other irreducible
component is isomorphic to $\cO\oplus \cO(-2)$.}
\end{figure}     

\subsection{Local Gromov--Witten invariants}
Let $Z \subset X$ be a closed subvariety of a smooth projective 
Calabi--Yau threefold.
Let $\Mbar_{g}(X,\vd)$ denote the stack of genus $g$ stable maps to $X$ representing 
$\vd \in H_{2}(X,\bZ)$. It is a Deligne--Mumford stack with a perfect obstruction
theory of virtual dimension zero which defines a virtual fundamental
zero-cycle $[\Mbar_{g} (X,\vd )]^{\vir}$. 

Whenever the substack $\Mbar_{g} (Z)$
consisting of stable maps whose image lies in $Z$ is a union of path connected
components of $\Mbar_{g} (X,\vd )$, it inherits a degree-zero virtual class. The
genus-$g$ \emph{local Gromov--Witten invariant} of $Z$ in $X$ is defined to be
the degree of this virtual class, and is denoted by $N_{\vd}^{g} (Z\subset X)$.
We write $N_{\vd}^{g} (Z)$ when the formal neighborhood is understood.

We will consider genus $g$, degree $\vd$ local Gromov--Witten invariants $N^g_\vd(Y^N)$,
where
$$
\vd=\sum_{j=1}^N \bigl(d_{1,j}[A_j]+ d_{2,j}[B_j]+ d_{3,j}[C_j]\bigr)\in H_2(Y^N;\bZ).
$$
For simplicity, we write $\vd=(\lo{\vd})$ where $\vd_i=(d_{i,1},\ldots,d_{i,N})$.
In this paper, we always assume $\vd$ is {\em effective} in the sense that $d_{i,j}\geq 0$.
We will show that the  local invariants $N^g_\vd(Y^N)$ are well defined 
in the following cases:
\begin{enumerate}
\item[(i)] (The minimal trivalent configuration)
$d_{1,1}=d_{2,1}=d_{3,1}=1$.
\item[(ii)] (A chain of rational curves)
$d_{1,1}>0$, $d_{2,j}=d_{3,j}=0$ for $1\leq j\leq N$. 
\end{enumerate} 
We will see in \fullref{sec:construction} that the formal neighborhood of $Y^N$ 
has a cyclic symmetry, so one can cyclically permute $\vd_1,\vd_2,\vd_3$ in 
Case (ii).
We show that in the above cases the local invariants $N^g_\vd(Y^N)$ are equal to certain
global or ordinary Gromov--Witten invariants of a blowup of $\bP^{3}$ at points
(\fullref{sec: local to global}), and we compute these ordinary invariants using 
the geometry of the Cremona transform (\fullref{sec:cremona}). 
To state our results, define constants $C_g$ by
\begin{equation}\label{eqn:Cg}
\sum_{g=0}^\infty C_g t^{2g}=
\left(\frac{t/2}{\sin(t/2)}\right)^2
=\sum_{g=0}^\infty\frac{|B_{2g}(2g-1)|}{(2g)!} t^{2g}.
\end{equation}

\begin{theo}[The minimal trivalent configuration]\label{thm:tri}
Suppose that 
$$
d_{1,1}=d_{2,1}=d_{3,1}=1.
$$
Then 
$$
N^g_{\vd}(Y^N)=\left\{ \begin{array}{ll}
C_g &  \text{if } 1=d_{i,1} \geq \cdots \geq d_{i,N} \geq 0
      \text{ for }\ i=1,2,3, \\
 0  & \text{otherwise.}   
\end{array}\right.
$$
\end{theo}

\begin{theo}[A chain of rational curves]\label{thm:chain}
Suppose that 
$$
\vd_1=(d_1,\ldots,d_N),\quad \vd_2=\vd_3=(0,\ldots,0),
$$
where $d_1>0$. Then
$$
N^g_{\vd}(Y^N)=\left\{ \begin{array}{ll}
C_g d^{2g-3} & 
\begin{array}{l}\text{if } d_1 = d_2 = \cdots = d_k = d > 0\text{ and}\\
 d_{k+1} = d_{k+2} = \cdots = d_N=0 \text{ for some }1\leq k\leq N
\end{array}\\[10pt]
 0  & \text{otherwise.}   
\end{array}\right.
$$
\end{theo}

Our results are new and add to the list of configurations of 
rational curves for which the local Gromov--Witten invariants are known.

The configuration in \fullref{thm:chain} is an $A_N$ curve. It is
interesting to compare \fullref{thm:chain} with the result
for a generic contractible $A_N$ curve $E=E_1\cup\cdots\cup E_N$
from Bryan--Katz--Leung \cite[Proposition 2.10]{BKaL}:

\begin{fact}[A generic contractible $A_N$ curve \cite{BKaL}] 
Assume $d_i>0$ for $i=1,\ldots,N$. Let $N_g(d_1,\ldots,d_N)$
denote genus $g$ local Gromov--Witten invariants of $E$ in the class
$\sum_{j=1}^N d_j[E_j]$. Then 
$$
N_g(d_1,\ldots,d_N)=\left\{ \begin{array}{ll}
C_g d^{2g-3}&  d_1=\cdots=d_N=d>0,\\
 0  & \textit{otherwise.}   
\end{array}\right.
$$
\end{fact}

Note that $Y^1$ is the closed topological vertex.
By the results in Faber--Pandharipande \cite{FP} and Bryan--Karp \cite{BK}, $N^g_{\lo{d}}(Y^1)$ is defined
in the following cases:
\begin{itemize}
\item[(iii)] (Super-rigid $\bP^1$) $d_1>0, d_2=d_3=0$ (and its cyclic permutation).
\item[(iv)] (The closed topological vertex) $\lo{d}>0$.
\end{itemize}

\begin{fact}[Super-rigid $\bP^1$ \cite{FP}] \label{thm:onePone} 
Suppose that $d>0$. Then
$$
N^g_{d,0,0}(Y^1)=N^g_{0,d,0}(Y^1)= N^g_{0,0,d}(Y^1)=C_g d^{2g-3}. 
$$
\end{fact}

\begin{fact}[The closed topological vertex \cite{BK}]  \label{thm:threePone}
Suppose that $\lo{d}>0$. Then
$$
N^g_{\lo{d}}(Y^1)=\left\{ \begin{array}{ll}
C_g d^{2g-3}&  d_1=d_2=d_3=d>0,\\
 0  & \textit{otherwise.}   
\end{array}\right.
$$
\end{fact}

\subsection{Formal Gromov--Witten invariants}

The minimal trivalent configuration $Y^N$ together with its formal neighborhood
is a nonsingular formal toric Calabi--Yau (FTCY) scheme $\hY^N$. The {\em formal Gromov--Witten invariants} 
$\tN^g_\vd(\hY^N)$ of $\hY^N$ are defined for all nonzero effective classes 
(see \fullref{sec:cases} and Bryan--Pandharipande \cite[Section 2.1]{BP04}).
Moreover, 
$$
\tN^g_\vd(\hY^N)=N^g_\vd(Y^N)
$$ 
in all the above cases (i)--(iv). Introduce formal variables $\lambda,t_{i,j}$
and define
$$
\tZ_N(\la;\vt)=\exp\Bigl(\sum_{g\geq 0} \sum_{\vd} \la^{2g-2} \tN_{g,\vd}(\hY^N) e^{-\vd\cdot \vt} \Bigr)
$$
where $\vd$ runs over all nonzero effective classes, and 
$$
\vt=(\lo{\vt}),\quad \vt_i=(t_{i,1},\ldots, t_{i,N}),\quad
\vd\cdot\vt=\sum_{i=1}^3\sum_{j=1}^N d_{i,j} t_{i,j}. 
$$
We call $\tZ_N(\la;\vt)$ the {\em partition function} of formal Gromov--Witten invariants of $\hY^N$.
It is the generating function of {\em disconnected} formal Gromov--Witten invariants of $\hY^N$.

In \fullref{sec:formal}, we will compute $\tZ_N(\la;\vt)$
by the mathematical theory of the topological vertex (see
Li--Liu--Liu--Zhou \cite{LLLZ}) and get the following
expression (\fullref{thm:tZ}):
\begin{multline}\label{eqn:tZ}
\tZ_N(\la;\vt)
=\exp\biggl(\sum_{n=1}^\infty \frac{1}{n[n]^2}\sum_{i=1}^3\sum_{2\leq k_1\leq k_2 \leq N}
e^{-n(t_{i,k_1}+\cdots +t_{i,k_2}) }\biggr) \\
\cdot\sum_{\vmu} \tilde{\cW}_{\vmu}(q) \prod_{i=1}^3
(-1)^{|\mu^i|} e^{-|\mu^i|t_{i,1}} s_{(\mu^i)^t}(u^i(q,\vt_i)).
\end{multline}
where $\vmu=(\up{\mu})$ is a triple of partitions, 
$q=e^{\sqrt{-1}\la}$, $[n]=q^{n/2}-q^{-n/2}$. The precise definitions of $\tilde{\cW}_{\vmu}(q)$
and $s_{(\mu^i)^t}(u^i(q,\vt_i))$ will be given in \fullref{sec:notation}. 
In particular, we will show that 
\begin{equation}\label{eqn:tZone}
\tZ_1(\la;\vt)= \sum_{\vmu}\tilde{\cW}_{\vmu}(q) 
\prod_{i=1}^3 (-1)^{|\mu^i|}e^{-|\mu^i|t_i}\cW_{(\mu^i)^t}(q)
=\exp\left(\sum_{n=1}^\infty\frac{Q_n(\vt)}{-n[n]^2}\right)
\end{equation}
where $\vt=(\lo{t})$, $\cW_\mu(q)$ is defined by 
\eqref{qdimension} in \fullref{sec:notation}, and 
\begin{multline}\label{eqn:Qn}
Q_n(\vt) = \\
e^{-nt_1}+e^{-nt_2}+e^{-nt_3}- e^{-n(t_1+t_2)}-e^{-n(t_2+t_3)}-e^{-n(t_3+t_1)}
+e^{-n(t_1+t_2+t_3)}.
\end{multline}
In \fullref{sec:vertex}, we will compute $\tZ_N(\la;\vt)$ by the 
physical theory of the topological vertex (see
Aganagic--Klemm--Mari\~no--Vafa \cite{AKMV})
and get the following expression (\fullref{thm:Z}):
\begin{multline}\label{eqn:Z}
Z_N(\la;\vt)=\exp\biggl(\sum_{n=1}^\infty \frac{1}{n[n]^2}\sum_{i=1}^3\sum_{2\leq k_1\leq k_2 \leq N}
e^{-n(t_{i,k_1}+\cdots +t_{i,k_2}) }\biggr) \\
\cdot \sum_{\vmu}\cW_{\vmu}(q) 
\prod_{i=1}^3 (-1)^{|\mu^i|} e^{-|\mu^i|t_{i,1}} s_{(\mu^i)^t}(u^i(q,\vt_i))
\end{multline}
where $\cW_\vmu(q)$ is defined by \eqref{topvertex} in \fullref{sec:notation}. 
In particular, we will show that (\fullref{thm:Zone}):
\begin{equation}\label{eqn:Zone}
Z_1(\la;\vt)= \sum_{\vmu}\cW_{\vmu}(q) \prod_{i=1}^3 (-1)^{|\mu^i|}e^{-|\mu^i|t_i}\cW_{(\mu^i)^t}(q)
=\exp\left(\sum_{n=1}^\infty\frac{Q_n(\vt)}{-n[n]^2}\right)
\end{equation}
The equivalence of the physical and mathematical theories of the topological vertex boils
down to the following combinatorial identity:
\begin{equation}\label{eqn:physics-math}
\cW_{\up{\mu}}(q)=\tilde{\cW}_{\up{\mu}}(q). 
\end{equation}
It is known that \eqref{eqn:physics-math} holds when one of the three partitions
is empty (see the work of Li, C-C M Liu, K Liu and Zhou
\cite{LLZ,LLLZ}). When none of the partitions is empty, Klemm has
checked all the cases where $|\mu^i|\leq 6$ by computer. Up to now, a 
mathematical proof of \eqref{eqn:physics-math} in full generality is not available. 
Equations \eqref{eqn:tZone} and \eqref{eqn:Zone} imply the following
result.

\begin{theo}\label{thm:WtW}
$$
 \sum_{\vmu} \cW_\vmu(q) \prod_{i=1}^3(-1)^{|\mu^i|}e^{-|\mu^i|t_i} \cW_{(\mu^i)^t}(q)
=\sum_{\vmu}\tilde{\cW}_\vmu(q) \prod_{i=1}^3 (-1)^{|\mu^i|}e^{-|\mu^i|t_i}\cW_{(\mu^i)^t}(q).
$$
\end{theo}

\fullref{thm:WtW} provides further evidence of \eqref{eqn:physics-math}
equating the vertex amplitudes derived from physical and mathematical theories of the topological
vertex.

\subsection{The topological vertex}\label{sec:notation}
In \cite{AKMV}, Aganagic, Klemm, Mari\~{n}o, and Vafa proposed
that Gromov--Witten invariants  of any toric Calabi--Yau threefold can be expressed in terms
of certain  relative invariants of its $\bC^3$ charts, called
{\em the topological vertex}. They suggested that
these local relative invariants should count holomorphic
maps from bordered Riemann surfaces to $\bC^3$ where the boundary
circles are mapped to three explicitly specified Lagrangian
submanifolds $\lo{L}$. The topological vertex depends on three 
partitions $\vmu=(\lo{\mu})$, where
$\mu^i$ corresponds to the winding numbers (the homology classes of boundary circles) in $L_i$. 
There is a symmetry on $\bC^3$ cyclically permuting $\lo{L}$, so one expects
the topological vertex to be symmetric under a cyclic permutation of the three partitions $\lo{\mu}$.

In \cite{AKMV} the topological vertex was computed by using the conjectural relation 
between open Gromov--Witten invariants on toric Calabi--Yau threefolds and
Chern--Simons invariants 
of knots and links. It has the following form:
\be \label{topvertex}
\cW_{\vmu}(q)=q^{\kappa_{\mu^2}/2 +\kappa_{\mu^3}/2}\sum_{\rho,\rho^1,\rho^3}
c^{\mu^1}_{\rho \rho^1}  c^{(\mu^3)^t}_{\rho \rho^3}
\frac{ \cW_{(\mu^2)^t \rho^1}(q) \cW_{\mu^2 \rho^3}(q)} {\cW_{\mu^2}(q)} .
\ee
In \eqref{topvertex}, $\mu^t$ denotes the partition transposed to $\mu$. 
The expression \eqref{topvertex} involves various quantities that we now define. $\kappa_{\mu}$ is given by
$$\kappa_{\mu}=\sum_i \mu_i(\mu_i -2 i+ 1).$$
The coefficients $c_{\mu \nu}^{\rho}$ are Littlewood--Richardson
coefficients. They can be defined in terms of Schur functions as follows
$$s_{\mu}s_{\nu}=\sum_{\rho}c_{\mu \nu}^{\rho} s_{\rho}.$$
Here, Schur functions are regarded as a basis for the ring $\Lambda$ of symmetric polynomials in an infinite 
number of variables. The quantity ${\cal W}_{\mu}(q)$ can be also defined in terms of Schur functions 
as follows:
\be
\label{qdimension}
{\mathcal W}_{\mu}(q)=s_{\mu}\bigl(x_i=q^{-i+\frac{1}{2}}\bigr).  
\ee
One can show that
\be
{\mathcal W}_{\mu^t}(q)=q^{-\kappa_{\mu}/2} {\mathcal W}_{\mu}(q).
\label{qdimtrans}
\ee
We also define, in analogy to skew Schur functions, 
\be
\label{skewqd}
{\cal W}_{\mu/\nu}(q)=\sum_{\lambda}c^{\mu}_{\nu \lambda} {\cal W}_{\lambda}(q).
\ee
Finally, ${\mathcal W}_{\mu \nu}(q)$ is defined by
\be \label{symw}
{\cal W}_{\mu \nu}(q)=q^{\kappa_{\mu}/2 + \kappa_{\nu}/2}\sum_{\lambda} {\cal W}_{\mu^t/\lambda}(q)
{\cal W}_{\nu^t/\lambda}(q).
\ee
This expression for $\cW_{\mu\nu}(q)$ is different from the one used originally in \cite{AKMV}. The 
fact that both agree follows from cyclicity of the vertex, and it has been
proved in detail by Zhou \cite{Z}. The expression for the vertex in terms of Schur functions 
is given in Okounkov--Reshetikhin--Vafa \cite{ORV} where the cyclicity
of the vertex is also proved.

In Li--Liu--Liu--Zhou \cite{LLLZ} the topological vertex was interpreted and defined as local relative invariants of 
a configuration $C_1\cup C_2\cup C_3$ of three $\bP^1$'s meeting at a point $p_0$
in a relative Calabi--Yau threefold $(Z,D_1,D_2,D_3)$,
where $K_Z+D_1+D_2+D_3 \cong \cO_Z$, $C_i$ intersects $D_i$ at a point $p_i \neq p_0$, 
and $C_i\cap D_j$ is empty
for $i\neq j$. The partition $\mu^i$ corresponds to the ramification pattern over $p_i$.
It is shown in \cite{LLLZ} that Gromov--Witten invariants of any toric Calabi--Yau
threefold (or more generally, formal Gromov--Witten invariants of formal toric
Calabi--Yau threefolds) can be expressed in terms of local relative invariants 
as described above, and the gluing rules coincide with those stated in
Aganagic--Klemm--Mari\~no--Vafa \cite{AKMV}.
The following expression of the vertex was derived in \cite{LLLZ}:
\begin{multline}\label{eqn:tW}
 \tilde{\cW}_{\vmu}(q)= q^{-(\kappa_{\mu^1}-2\kappa_{\mu^2}-\frac{1}{2}\kappa_{\mu^3})/2}
\!\!\!\!\!\!\!\!\sum_{\nu^+,\nu^1,\nu^3,\eta^1,\eta^3 } c^{\nu^+}_{(\nu^1)^t\mu^2}
c_{(\eta^1)^t\nu^1}^{\mu^1}c_{\eta^3(\nu^3)^t}^{\mu^3}\\
\cdot q^{(-2\kappa_{\nu^+}  - \frac{\kappa_{\nu^3}}{2})/2} \cW_{\nu^+ \nu^3}(q)
\sum_{\si} \frac{1}{z_\si}\chi_{\eta^1}(\sigma)\chi_{\eta^3}(2\sigma).
\end{multline}
Here $2\si=(2\si_1\geq 2\si_2\geq \cdots)$ if $\si=(\si_1\geq\si_2\geq \cdots)$.
Recall that 
$$
z_\sigma=\prod_{i\geq 1}{i^{m_i}\cdot m_i!}
$$ where $m_i=m_i(\sigma)$ is the number of parts of the partition
$\sigma$ equal to $i$ (see Macdonald \cite[p.17]{M}).

It is expected that the two different enumerative interpretations
in \cite{AKMV} and in \cite{LLLZ} of the vertex
give rise to equivalent counting problems, in the spirit of the following simple
example: counting ramified covers of a disc by bordered Riemann surfaces with 
prescribed winding numbers is equivalent to 
counting ramified covers of a sphere by closed Riemann surfaces with prescribed 
ramification pattern over $\infty$.

Finally, we introduce some notation which will arise in computations in
\fullref{sec:vertex}.
For any positive integer $n$, define
\begin{equation}\label{eqn:un}
u^i_n(q,\vt_i)=\frac{1}{[n]}\Bigl(1+\sum_{k=2}^N e^{-n(t_{i,2}+\cdots+ t_{i,k})}\Bigr).
\end{equation}
Given a partition 
$\mu=(\mu_1\geq \mu_2 \geq \cdots\geq \mu_{\ell}>0 )$, define
\begin{equation}\label{eqn:umu}
u^i_\mu(q,\vt_i)=\prod_{j=1}^{\ell} u_{\mu_j}^i(q,\vt_i).
\end{equation}
and
\begin{equation}\label{eqn:schur-u}
s_\mu(u^i(q,\vt_i))=\sum_{|\nu|=|\mu|} \frac{\chi_\mu(\nu)}{z_\nu} u^i_\nu(q,\vt_i).
\end{equation}
In particular, when $N=1$, we have
\[
u^i_n=\frac{1}{[n]}=\sum_{i>0} q^{-i+1/2}
\]
So
\begin{equation}\label{eqn:schur-W}
s_\mu(u^i(q,\vt_i))=s_\mu(x_i=q^{-i+\frac{1}{2} })=\cW_\mu(q).
\end{equation}

\subsection*{Acknowledgements}
The authors give warm thanks to Jim Bryan for helpful conversations. 
The first author thanks NSERC for its support.
The second author thanks Jun Li, Kefeng Liu, Jian Zhou for collaboration
\cite{LLLZ} and Shing-Tung Yau for encouragement. 
Finally, the authors thank the referee for the detailed comments and
numerous valuable suggestions on the presentation.
\section{Cremona}\label{sec:cremona}
In this section we prove Theorems \ref{thm:tri} and \ref{thm:chain} using the
geometry of the Cremona transform. We assume that the formal neighborhood 
$Y^{N}\subset X$ is as constructed in \fullref{sec:construction}. We also
assume that the local invariants of $Y^{N}$ are equal to certain ordinary
invariants of $X$, which we prove in \fullref{sec: local to global}.

\subsection{The blowup of $\bC\bP^{3}$ at points}\label{subsection: the blowup of cp3 at points}

We briefly review the properties of the blowup of $\bP^{3}$ at points used here
for completeness and to set notation. This material can be found in much 
greater detail in, for instance, Griffiths--Harris \cite{GH}.

Let $X\rightarrow \bP^{3}$ be the blowup of $\bP^{3}$ along $M$ distinct
points $\{p_{1},\ldots ,p_{M}\}$. We describe the homology of $X$.
All (co)homology is taken with integer coefficients.
Note that we may identify homology and cohomology 
as rings via Poincar\'e duality, where cup product is dual to intersection
product.

Let $H$ be the total
transform of a hyperplane in $\bP^{3}$, and let $E_{i}$ be the exceptional
divisor over $p_{i}$. Then $H_{4} (X,\bZ)$  has a basis
\[
H_{4} (X) = \left<H,E_{1},\ldots ,E_{M}    \right>.
\]
Furthermore, let $h \in H_{2} (X)$ be the class of a line in $H$, and let
$e_{i}$ be the class of a line in $E_{i}$. The collection of all such classes
form a basis of $H_{2} (X)$.
\[
H_{2} (X)= \left<h,e_{1},\ldots ,e_{M} \right>
\]
The intersection ring structure is given as follows.
Let $pt \in H_{0} (X) $ denote the
class of a point. Two general hyperplanes meet in a line,
so $H\cdot H= h$. A general hyperplane and line intersect in a point, so
$H \cdot h=pt$. Also, a general hyperplane is far from the center
of a blowup, so all other products involving $H$ or $h$ vanish. The
restriction of $\cO_{X} (E_{i})$ to $E_{i}\cong \bP^{2}$ is the dual of
the bundle $\cO_{\bP^{2}} (1)$, so $E_{i}\cdot E_{i}$ is represented
by minus a hyperplane in $E_{i}$, i.e. 
$E_{i}\cdot E_{i} =-e_{i}$, and 
$E_{i}^{3}= (-1)^{3-1}pt=pt$ (see Fulton \cite{Fu}). Furthermore, the
centers of the
blowups are far away from each other, so all other intersections vanish.
In summary, the following are the only non-zero intersection products.
\begin{equation*}
\boxed{\begin{aligned}
	H \cdot H &= h  &H\cdot h &=pt\\ 
  E_{i}\cdot E_{i} &= -e_{i} &E_{i}\cdot e_{i} &=-pt\\
\end{aligned}}
\end{equation*}
Also, we point out the that the canonical bundle $K_{X}$ is easy to
describe in this basis:
\[
K_{X} = -4H +2\sum_{i=1}^{M}E_{i}
\]
Finally, we introduce a notational convenience for the Gromov--Witten invariants
of $\bP^{3}$ blown up at points in a Calabi--Yau class. Any curve class is of 
the form
\[
\beta = dh - \sum_{i=1}^{M}a_{i}e_{i} 
\]
for some integers $d,a_{i}$ where $d$ is non-negative.
Thus $K_{X}\cdot \beta =0$ if and only
if $2d=\sum_{i=1}^{M} a_{i}$. In that case,
the virtual dimension of 
$\Mbar_g(X,\beta)$ is zero, and 
$$
\left< \; \right>_{g,\beta}^{X}=\int_{[\Mbar_g(X,\beta)]^{\vir} } 1
$$ is determined by the discrete data
$\{d,a_{i},\ldots ,a_{M} \}$.
Then, we may use the shorthand notation 
\[
\left< \; \right>_{g,\beta}^{X} = \left<d; a_{1},\ldots ,a_{M} \right>^{X}_{g}.
\]
For example,
\[
\left< \; \right>_{g,5h -e_{1}-e_{2} -2e_{3}-3e_{5}-3e_{6} }^{X} =
\left< 5; 1,1,2,0,3,3 \right>_{g}^{X}.
\]
Furthermore, the Gromov--Witten invariants of $X$ do not depend on ordering of 
the points $p_{i}$, and thus for any permutation $\sigma$ of $M$ points,
\[
\left<d; a_{1},\ldots, a_{M} \right>^{X}_{g} =
\left<d; a_{\sigma (1)},\ldots, a_{\sigma (M)} \right> ^{X}_{g}. 
\]

\subsection{Properties of the invariants of the blowup of $\bP^{3}$ at points}
\label{subsec: properties of P3}
First, we use the fact, shown in Bryan--Karp \cite{BK}, that the 
Gromov--Witten invariants of
the blowup of $\bP^{3}$ along points have a symmetry which arises from
the geometry of the Cremona transformation.

\begin{theo}[Bryan--Karp \cite{BK}]\label{thm: cremona invariance}
Let $ \beta =dh -\sum _{i=1}^{M}a_{i}e_{i}$ with $2d=\sum _{i=1}^{M}a_{i}$
and assume that $a_{i}\neq 0$ for some $i>4$. Then we have the following
equality of Gromov--Witten invariants:
\[
\left\langle \; \right\rangle^{X}_{g,\beta }=\left\langle \;
\right\rangle^{X}_{g,\beta '}
\]
where $\beta '=d'h-\sum _{i=1}^{M}a_{i}'e_{i}$  has coefficients given by
\begin{align*}
d'\, &=3d-2 (a_{1}+a_{2}+a_{3}+a_{4})\\
a_{1}'&=\;\, d-\; \, (a_{2}+a_{3}+a_{4})\\
a_{2}'&=\;\, d- \;\, (a_{1}+a_{3}+a_{4})\\
a_{3}'&=\;\, d-\; \, (a_{1}+a_{2}+a_{4})\\
a_{4}'&=\;\, d-\; \, (a_{1}+a_{2}+a_{3})\\
a_{5}'&=\;\, a_{5}\\
&\;\; \vdots &\\
a_{M}'&=\;\, a_{M}.
\end{align*}
\end{theo}

We also use the following vanishing lemma, and a few of its corollaries. 

\begin{lemm}\label{lem: vanishing}
Let $X$ be the blowup of $\bP^3$ at $M$ distinct generic points 
$\{ x_1, \ldots, x_M\}$,
and $\beta=dh-\sum_{i=1}^{M} a_i e_i$ with 
$2d=\sum_{i=1}^{M} a_i$, and assume that $d>0$ and $a_i <0$ for some $i$. Then
\[
\Mbar_g(X,\beta) = \emptyset.
\]
\end{lemm}

\begin{coro}\label{cor: vanishing}
For any $M$ points $\{ x_1, \ldots, x_{M}\}$ and $X$ and $\beta$ as 
above the
corresponding invariant vanishes;
\[
\left\langle \; \right\rangle^{X}_{g,\beta } =0.
\]
\end{coro}

This follows immediately from the deformation invariance of Gromov--Witten 
invariants and \fullref{lem: vanishing}.

\begin{proof} 
In genus zero, \fullref{lem: vanishing} follows from a vanishing theorem of
Gathmann \cite[Section~3]{G}.
In order to prove \fullref{lem: vanishing}, for arbitrary genus, 
it suffices to
show that the result holds for a specific choice of points, as if the
moduli space is empty for a specific choice, then it is empty for the generic
choice. By choosing some of the points to be coplanar, and the rest to also be 
coplanar on a second plane, the result follows. For further details, see
Karp \cite{K}.
\end{proof}

\begin{coro}\label{cor: extraneous blowups}
Let $X$ be the blowup of $\bP^{3}$ along $M$ points and define
$\beta =dh -\sum_{i=1}^{M}a_{i}e_{i}$ where $2d=\sum_{i=1}^{M}a_{i}$
and $d>0$.
Also define
\[
X' \xrightarrow{\pi} X
\]
to be the blowup of $X$ at a generic point $p$, so that $X'$ is deformation 
equivalent to the blowup of $\bP^{3}$ at $M+1$ distinct points.
Let $\{h',e'_{1},\ldots ,e'_{M+1} \}$ be a basis of $H_{2} (X')$, and let
$\beta' = dh'-\sum_{i=1}^{M}a_{i}e'_{i}$.
Then
\[
\left<d;a_{1},\ldots a_{M},0 \right>_{g}^{X'} 
=\left<d;a_{1},\ldots ,a_{M} \right>_{g}^{X}
\]
\end{coro}

\begin{proof}
This result follows from the more general results of Hu \cite{H}. An 
independent proof using \fullref{lem: vanishing} can be found in
Karp \cite{K}.
\end{proof}

\subsection[Proof of \ref{thm:tri}]
  {Proof of \fullref{thm:tri}}\label{subsec: proof tri}
Let the blowup space $X^{N+1}$ and
the minimal trivalent configuration $Y^{N}$ be as constructed in
\fullref{sec:construction} on page~\pageref{sec:construction}. By 
\fullref{thm:triIm} on page~\pageref{thm:triIm}  we have 
\[
N^{g}_{\vd} (Y^{N}) =\braket^{X^{N+1}}_{g,\vd}.
\]
Assume that the invariant is non-zero:
\begin{align*}
\braket_{g,\vd}^{X^{N+1}}=&\bigl<3;1,1-d_{1,2},\dotsc ,d_{1,N-1}-d_{1,N},d_{1,N},\\
                &\;  1,1-d_{2,2},\dotsc ,d_{2,N-1}-d_{2,N},d_{2,N},\\
                &\;  1,1-d_{3,2},\dotsc ,d_{3,N-1}-d_{3,N},d_{3,N} \bigr>^{X^{N+1}}_g\\
                & \;  \neq 0
\end{align*}
Then, by \fullref{cor: vanishing} , the coefficient of each 
$e_{i},f_{i},g_{i}$ is non-negative. Thus, for $i=1,2,3$,
\begin{equation}\label{eqn:tri-nonvanishing}
1 \geq d_{i,2} \geq \dotsb \geq d_{i,N} \geq 0.
\end{equation}
Therefore we compute
\begin{align*}
\braket_{g,\vd}^{X^{N+1}} =& \bigl< 3; 1,0,\dotsc ,0,1,\\
                           & \; 1,0,\dotsc ,0,1,\\
                           & \; 1,0,\dotsc ,0,1\bigr>_{g}^{X^{N+1}}\\
                           &= \left<3;1,1,1,1,1,1 \right>^{X^{2}}_{g},
\end{align*}
where the last equality follows from \fullref{cor: extraneous blowups}.
So when \eqref{eqn:tri-nonvanishing} holds, we have
$$
N^g_\vd(Y^N)=N^g_{1,1,1}(Y^1)=C_g.
$$
The last equality follows from \fullref{thm:threePone} (see Bryan--Karp
\cite{BK}).
\qed

\subsection[Proof of \ref{thm:chain}]
  {Proof of \fullref{thm:chain}}\label{subsec: proof chain}
Let the blowup space $\tX^{N+1}$ and the chain of rational curves $Y^{N}_A$ be
as constructed in \fullref{sec:construction} on page~\pageref{eq:chain}. 
By \fullref{thm:chainIm} on page~\pageref{thm:chainIm}  we have
\[
N^g_{\vd}(Y^N)=N^g_{\vd_1}(Y^N_A)  =\braket^{\tX^{N+1}}_{g,\vd}
\]
where 
$$\vd_1=(d_1,\ldots,d_N), \quad \vd_2=\vd_3=(0,\ldots,0).$$
Assume that the invariant is non-zero:
\[
\braket^{\tX^{N+1}}_{g,\vd}
=\left<d_{1};d_{1},d_{1}-d_{2},\dotsc ,d_{N-1}-d_{N},d_{N} \right>_g^{\tX^{N+1}}
\neq 0.
\]
By \fullref{cor: vanishing} the multiplicities are decreasing:
\[
d_{1}\geq d_{2} \geq \dotsb \geq d_{N} \geq 0
\]
Therefore, as $d_{1}>0$, there exists some $1\leq j\leq N$ such that 
\[
d_{1}\geq \dotsb \geq d_{j}>0, \quad  d_{j+1}=\dotsb =d_{N}=0.
\]
Then, using \fullref{cor: extraneous blowups}, we compute
\begin{align*}
N^{g}_{\vd} (Y^{N})&= \left<d_{1};d_{1},d_{1}-d_{2},\dotsc ,d_{j-1}-d_{j},0,\dotsc 0 \right>_g^{\tX^{N+1}}\\
 &=\left<d_{1};d_{1},d_{1}-d_{2},\dotsc ,d_{j-1}-d_{j},d_{j} \right>^{\tX^{j+1}}_{g}.
\end{align*}
Note that for any $1\leq i \leq j+1$ we may reorder
\begin{multline*}
\left<d_{1};d_{1},d_{1}-d_{2},\dotsc ,d_{j-1}-d_{j},d_{j} \right>^{\tX^{j+1}}_{g}=\\
\left<d_{1};d_{1},d_{1},d_{i}-d_{i+1},0,0,d_{1}-d_{2},\dotsc\right.\\
\left.\dotsc,d_{i-2}-d_{i-1},
d_{i+1}-d_{i+2},\dotsc ,d_{j-1}-d_{j},d_{j} \right>^{\tX^{j+1}}_{g}
\end{multline*}
Applying Cremona invariance (\fullref{thm: cremona invariance}) we
compute
\begin{multline*}
\braket^{\tX^{N+1}}_{g,\vd} = \bigl< d_{1}-2 (d_{i}-d_{i+1}); d_{1}-
(d_{i}-d_{i+1}),0,d_{i+1}-d_{i}, d_{i+1}-d_{i},\\
d_{1}-d_{2},\dotsc ,d_{j-1}-d_{j},d_{j}\bigr>^{\tX^{j+3}}_{g}.
\end{multline*}
Then, by \fullref{cor: vanishing}, $d_{i+1}\geq d_{i}$. Since this
inequality holds for every $1\leq i \leq j$ we have 
$d_{1}\leq \dotsb \leq d_{j}$. Therefore
\[
d_{1}=\dotsb =d_{j}=d.
\]
Thus we have
\begin{align*}
\braket^{\tX^{N+1}}_{g,\vd}&= \left<d;d,0,\dotsc ,0,d \right>^{\tX^{j+1}}_{g}\\
                           &=\left<d;d,d \right>^{\tX^{2}}_{g}\\
                           &= N^g_{d,0,0}(Y^1)\\
                           &= C_g d^{2g-3}
\end{align*}
The last equality follows from Faber--Pandharipande \cite{FP}.
\qed 

\section{Construction}\label{sec:construction}

We construct these configurations as 
subvarieties of a locally Calabi--Yau space $X^{N+1}$, which is obtained via a
sequence of toric blowups of $\bP^{3}$:
\[
 X^{N+1} \xrightarrow{\pi_{N+1}} X^{N} \xrightarrow{\pi_{N}}
\cdots \xrightarrow{\pi_{2}} X^{1} \xrightarrow{\pi_{1}} X^{0} =\bP^{3}
\]
In fact, $X^{i+1}$ will be the blowup of
$X^{i}$ along three points.
Our rational curves will be labeled by $A_{i}, B_{i},C_{i}$,
where $1\leq i \leq N$, reflecting the nature of
the configuration. Curves in intermediary spaces will 
have super-scripts, and their corresponding proper transforms in $X$ will not.

The standard torus $\bT = (\bC^{\times})^{3}$ action on $\bP^{3}$
is given by
\[
(t_{1}, t_{2} , t_{3}) \cdot (x_{0}\colon x_{1}\colon x_{2} \colon x_{3}) \mapsto
(x_{0}\colon t_{1}x_{1}\colon t_{2}x_{2} \colon t_{3}x_{3}).
\]
There are four $\bT$--fixed points in $X^{0}\co=\bP^{3}$; we label them
$p_{0}= (1\colon 0\colon 0\colon 0)$, $q_{0}=(0 \colon  1 \colon  0 \colon  0)$,
$r_{0}= (0\colon 0\colon 1 \colon 0)$ and 
$s_{0}= (0\colon 0\colon 0\colon 1)$. Let $A^{0}$, $B^{0}$ and $C^{0}$
denote the (unique, $\bT$--invariant) line in $X^{0}$ through the two points
 $\{p_{0},s_{0}\}$, $\{q_{0},s_{0}\}$ and $\{r_{0},s_{0}\}$, respectively. 

Define
\[ 
X^{1} \xrightarrow{\pi_{1}} X^{0}
\]
to be the blowup of $X^{0}$ at the three 
points $\{p_{0},q_{0},r_{0}\}$,  and let $A^{1},B^{1},C^{1} \subset X^{1}$
be the proper 
transforms of $A^{0},B^{0}$ and $C^{0}$. The exceptional divisor in
$ X^{1}$ over $p_{0}$
intersects $A^{1}$ in a unique fixed point; call it $p_{1} \in X^{1}$. 
Similarly, the
exceptional divisor in $X^{1}$ also intersects each of
$B^{1}$ and $C^{1}$ in unique fixed points; call them  $q_{1}$ and $r_{1}$.

\begin{figure}[ht]
\label{fig: X2}
\centerline{\includegraphics{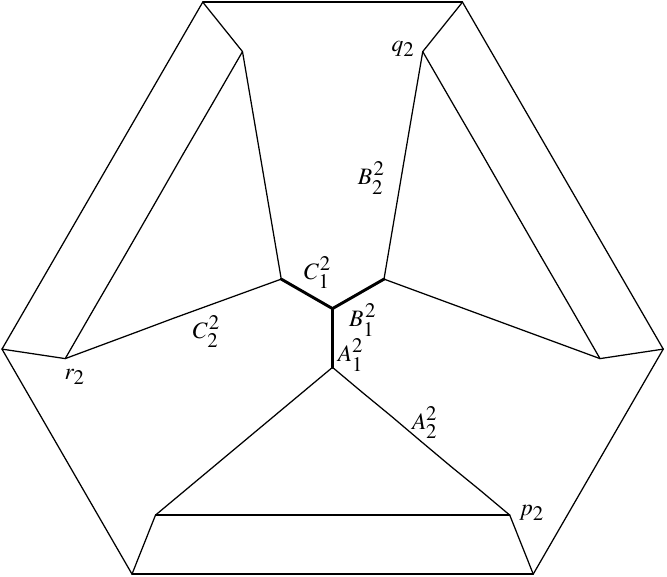}}
\caption{The $\bT$--invariant curves in $X^{2}$} 
\end{figure}
Now define 
\[
X^{2} \xrightarrow{\pi_{2}} X^{1} 
\]
to be the blowup
of $X^{1}$ at the three points $\{p_{1},q_{1},r_{1}\}$, and let 
$A_{1}^{2},B_{1}^{2},C_{1}^{2} \subset X^{2}$ be the proper
transforms of $A^{1},B^{1},C^{1}$. The exceptional divisor over $p_{1}$ 
contains two
$\bT$--fixed points disjoint from $A^{2}_{1}$. Choose one of them, and call it 
$p_{2}$; this choice is arbitrary. Similarly, there are two fixed points in the
exceptional divisors above $q_{1},r_{1}$ disjoint from $B_{1}^{2},C^{2}_{1}$. 
Choose
one in each pair identical to the choice of $p_{2}$ and call them $q_{2}$ 
and $r_{2}$ (identical makes sense here as the configuration of curves in 
\fullref{fig: X2} is rotationally symmetric). 
This choice is indicated in \fullref{fig: X2}. Let
$A^{2}_{2}$  denote the (unique, $\bT$--invariant) line
intersecting $A^{2}_{1}$ and $p_{2}$. Define $B^{2}_{2}, C^{2}_{2}$ 
analogously.

Clearly $X^{2}$ is deformation equivalent to a 
blowup of 
$\bP^{3}$ at six distinct points. The $\bT$--invariant curves in $X^{2}$ are 
depicted  in
\fullref{fig: X2}, where each edge corresponds to a 
$\bT$--invariant 
curve in $X^{2}$, and each vertex corresponds to a fixed point.

\begin{figure}[h]
\centerline{\includegraphics{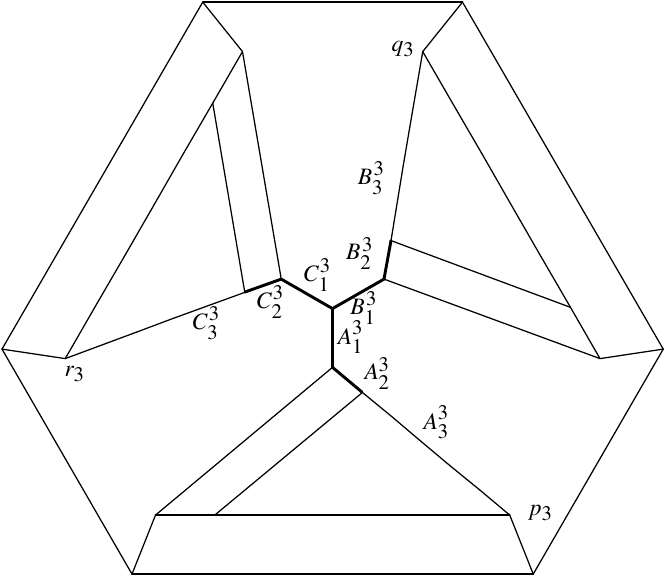}}
\label{fig: X3}
\caption{The $\bT$--invariant curves in $X^{3}$}
\end{figure}

We now define a sequence of blowups beginning with $X^{2}$. Fix an 
integer $N \geq 2$. For each $1 < i \leq N$, define
\[
X^{i+1} \xrightarrow{\pi_{i+1}} X^{i}
\]
to be the blowup of $X^{i}$ along the three points $p_{i},q_{i},r_{i}$.
Let $A^{i+1}_{j}\subset X^{i+1} $ denote the proper transform of $A^{i}_{j}$ 
for 
each $1 \leq j \leq i $. The
exceptional divisor in $X^{i+1}$ above $p_{i}$ contains two $\bT$--fixed points,
choose one of them and call it $p_{i+1}$. Similarly choose $q_{i+1}, r_{i+1}$,
and define $A^{i+1}_{i+1} \subset X^{i+1}$ to be the line intersecting
$A^{i+1}_{i}$ and $p_{i+1}$, with $B^{i+1}_{i+1}, C^{i+1}_{i+1}$ defined
similarly. The $\bT$--invariant curves in $X^{3}$ are shown in 
\fullref{fig: X3}. 

\begin{figure}[ht!]
\label{fig: XN+1}
\centerline{\includegraphics{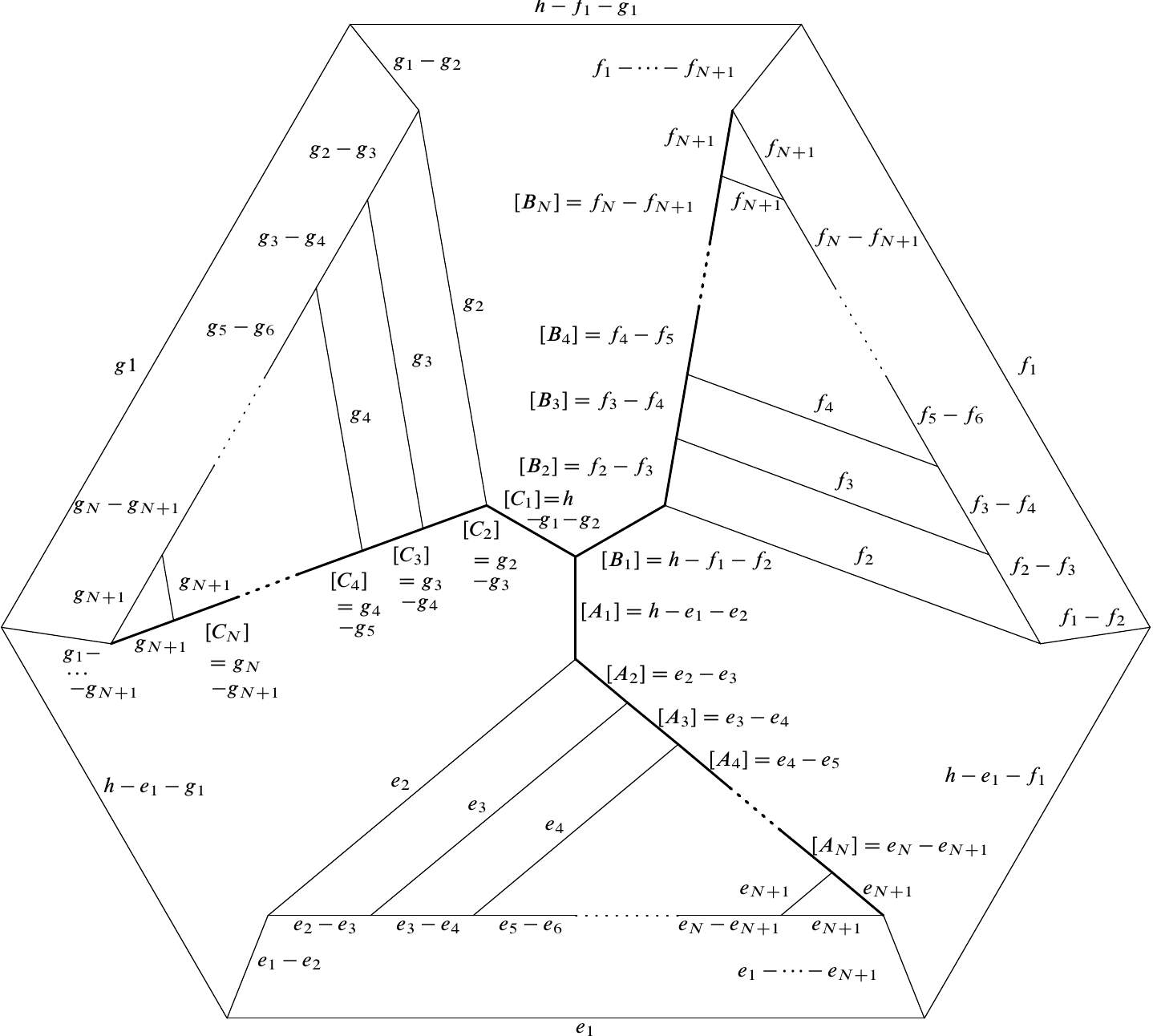}}
\caption{The $\bT$--invariant curves in $X^{N+1}$}
\end{figure}

Finally, we define the {\em minimal trivalent configuration} $Y^{N} \subset X^{N+1}$ by
\[
Y^{N} = \bigcup_{1 \leq j \leq N } A_{j} \cup B_{j}\cup 
C_{j},
\]
where
\[
A_j=A^{N+1}_j,\quad B_j=B^{N+1}_j,\quad C_j=C^{N+1}_j.
\]
The configuration $Y^{N}$ is shown in \fullref{fig: XN+1}, along with 
all other $\bT$--invariant curves in $X^{N+1}$. It contains a chain
of rational curves:
\begin{equation*}\label{eq:chain}
Y^N_A=A_1\cup\cdots\cup A_N.
\end{equation*}

\subsection{Homology}\label{subsec: homology} 
We now compute $H_{*} (X^{N+1},\bZ)$ and identify the class of the 
configuration
$[Y^{N}]\in H_{2} (X^{N+1},\bZ)$. All (co)homology will be taken with
 integer 
coefficients. We denote
divisors by upper case letters, and curve classes with the lower case. In 
addition, we decorate homology classes in intermediary spaces with a tilde,
and their total transforms in $X^{N}$ are undecorated.

Let $\tilde{E}_{1},\tilde{F}_{1},\tilde{G}_{1} \in H_{4} (X^{1})$
denote the exceptional divisors in $X^{1} \rightarrow X^{0}$ over the
points $p_{0},q_{0}$ and $r_{0}$, and let $E_{1},F_{1},G_{1} \in H_{4}
(X)$ denote their total transforms.  Continuing, for each $1 \leq i
\leq N+1$, let $\tilde{E}_{i},\tilde{F}_{i},\tilde{G}_{i} \in H_{4}
(X^{i})$ denote the exceptional divisors over the points
$p_{i-1},q_{i-1},r_{i-1}$ and let $E_{i},F_{i},G_{i} \in H_{4} (X) $
denote their total transforms. Finally, let $H$ denote the total transform of
the hyperplane in $X^{0}=\bP^{3}$.
The collection of all such classes
$\{H,E_{i},F_{i},G_{i} \}$, where $1 \leq i \leq N+1$, spans $H_{4}
(X^{N+1})$.

Similarly, for each $1 \leq i \leq N+1$, let
$
\tilde{e}_{i},\tilde{f}_{i},\tilde{g}_{i} \in H_{2} (X^{i+1})$ 
denote the class of a line in 
$\tilde{E}_{i},\tilde{F}_{i},\tilde{G}_{i}$
and let 
$e_{i},f_{i},g_{i} \in H_{2} (X)$ denote their total transforms.
In addition,  let $h\in H_{2} ( X^{N+1})$ denote the class of a line in $H$.
Then $H_{2} (X^{N+1})$ has a basis given by
$\{h,e_{i},f_{i},g_{i}\}$.

The intersection product
ring structure is given as follows. Note that $X^{N+1}$ is deformation
equivalent to the blowup of $\bP^{3}$ at 
$3N$ distinct points. Therefore, these
\begin{equation*}
\boxed{\begin{aligned}
	H \cdot H &= h  &H\cdot h &=pt\\ 
  E_{i}\cdot E_{i} &= -e_{i} &E_{i}\cdot e_{i} &=-pt\\
  F_{i}\cdot F_{i} &= -f_{i} &F_{i}\cdot f_{i} &=-pt\\
  G_{i}\cdot G_{i} &= -g_{i} &G_{i}\cdot g_{i} &=-pt
\end{aligned}}
\end{equation*}
are all of the nonzero intersection products in $H_{*} (X^{N+1})$.

In this basis, the classes of the components of $Y^{N}$ are given as follows.
\begin{align*}
&[A_{i}]=
\begin{cases}
h-e_{1}-e_{2} & \text{if } i=1\\
e_{i}-e_{i+1} & \text{otherwise}  
\end{cases}\\
&[B_{i}]=
\begin{cases}
h-f_{1}-f_{2} & \text{if } i=1\\
f_{i}-f_{i+1} & \text{otherwise}  
\end{cases}\\
&[C_{i}]=
\begin{cases}
h-g_{1}-g_{2} & \text{if } i=1\\
g_{i}-g_{i+1} & \text{otherwise}  
\end{cases}
\end{align*}
To see this, recall that $A_{1}$ is the proper transform of a line through two
points which are centers of a blowup, and that $A_{i}$, for $i>1$, is the 
proper transform of a line in an exceptional divisor containing a center of a
blowup. $B_{i}$ and $C_{i}$ are similar.

\section{Local to global}\label{sec: local to global}
In this section, we will show that 
the local invariants $N^{g}_{\vd}(Y)$ are equal to the ordinary
invariants $\braket_{g,\vd}^{X^{N+1}}$ in case $Y$ is either the
minimal trivalent configuration $Y^N$ or  
the  chain of rational curves $Y_A^N$ defined in \fullref{sec:construction}.

\subsection{The minimal trivalent configuration}\label{subsec: local to global mtc}

\begin{prop}\label{thm:triIm}
Let $f\co  \Sigma \to X^{N+1}$ represent a point in 
$\Mbar_g(X^{N+1},\vd)$, where
$$
1=d_{i,1}\geq \cdots \geq d_{i,N} \geq 0. 
$$
Then the image of $f$ is contained in the minimal
trivalent configuration 
\[
Y^N=\bigcup_{1\leq j\leq N} A_j\cup B_j\cup  C_j.
\]
\end{prop}

\begin{proof} 
We use the toric nature of 
the construction. Assume that there exists a stable map
\[
[f\co \Sigma \rightarrow X^{N+1}] \in \Mbar_g(X^{N+1},\vd)
\]
such that $\im(f) \not\subset Y^{N}$. Then there exists a point 
$p\in \im(f)$ such that $p \not\in Y^{N}$.

Recall that $\bT$--invariant
subvarieties of a toric variety are given precisely by orbit closures of
one-parameter subgroups of $\bT$. So in particular the limit of $p$ under
the action of a one-parameter subgroup is a $\bT$--fixed point. Moreover,
since  $p \not\in Y^{N}$, there exists a one-parameter subgroup
$\psi\co  \bC^{\times} \rightarrow \bT$ such that
\[
\lim_{t\rightarrow 0} \psi (t) \cdot p = q
\]
where $q$ is $\bT$--fixed and $q \not\in Y^{N}$.

The limit of $\psi$ acting on $[f]$ is a stable map $f'$ such that 
$q \in \im(f' )$. It follows that $q$ is in the image of all stable maps in
the orbit closure of $[f']$. Thus, there must exist a stable map 
$[f''\co  \Sigma \rightarrow X^{N+1}] \in \Mbar_g(X^{N+1},\vd)$ such that
$\im(f'')$ is $\bT$--invariant and $\im(f'')\not\in Y^{N}$.

We show that this leads to a contradiction. Let $F$ denote the union of the
$\bT$--invariant curves in $X^{N+1}$; it is shown above in 
\fullref{fig: XN+1}. We study the possible components of $F$ contained in
the image of $f''$.

\begin{figure}[ht!]
\centerline{\includegraphics{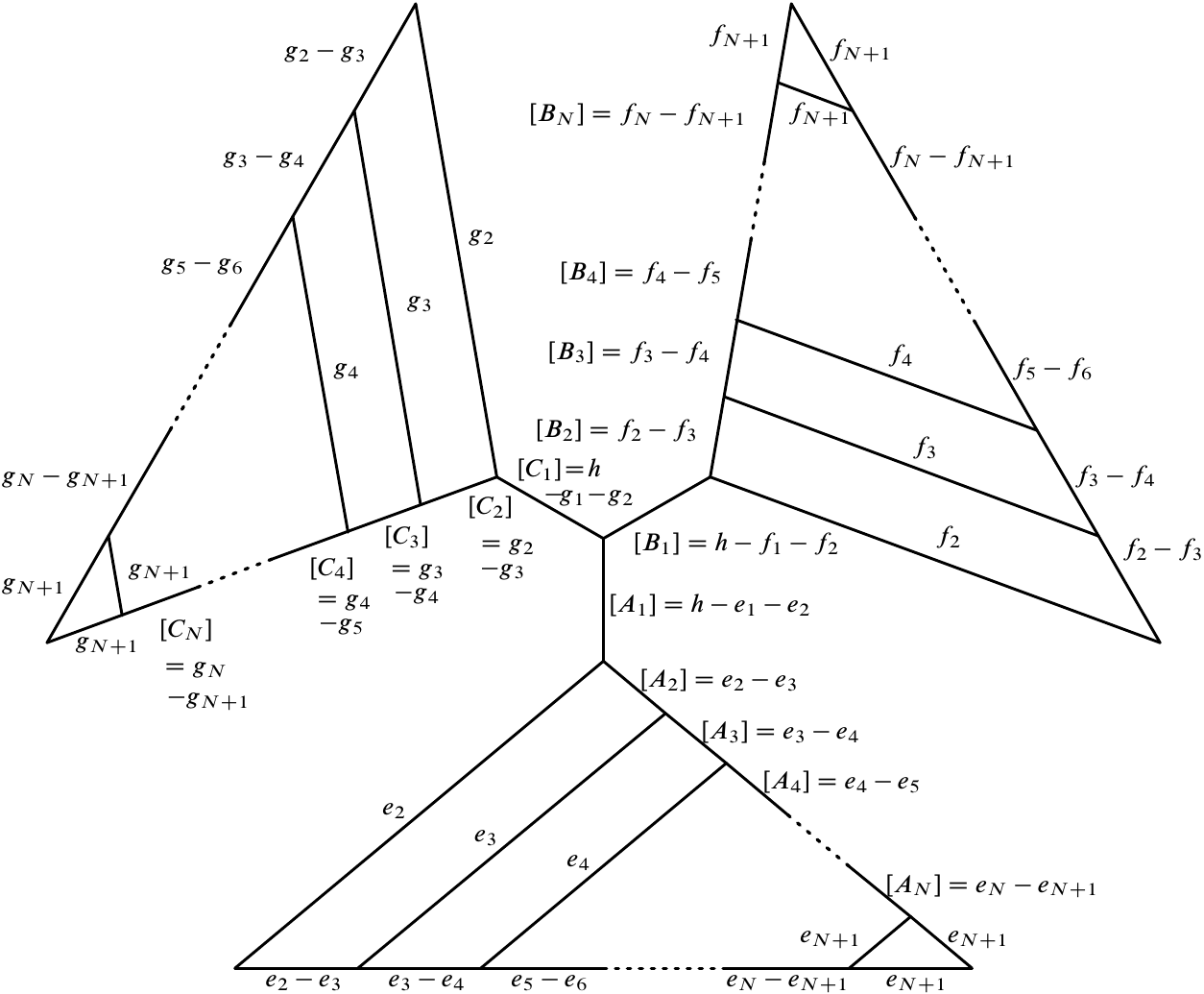}}
\caption{The possible curves in $\im(f'')$}
\label{fig: no h's or e's}
\end{figure}

Note that the push forward of the class of $\Sigma$ is given by
\begin{align*}
f''_{*} [\Sigma]  = \, 3h 
    &  - e_{1} - \sum_{j=1 }^{N-1} (d_{1,j}-d_{1,j+1})e_{j+1} -d_{1,N}e_{N+1}\\
 & - f_{1} - \sum_{j=1 }^{N-1} (d_{2,j}-d_{2,j+1})f_{j+1} -d_{2,N}f_{N+1}\\
 & - g_{1} - \sum_{j=1 }^{N-1} (d_{3,j}-d_{3,j+1})g_{j+1} - d_{3,N}g_{N+1}.
\end{align*}
Suppose that $A_{1}\cup B_{1}\cup C_{1} \subset \im(f'')$. Then 
$f''_{*} [\Sigma]$ contains (at least) $3h$. Note that $[F]$ has no $-h$ terms.
Therefore
$\im(f'')$ does not contain any of the curves 
$h-e_{1}-f_{1}, h-e_{1}-g_{1}, h-f_{1}-g_{1}$. And furthermore each of $A_{1}$,
$B_{1}$ and $C_{1}$ must have multiplicity one.

There are no remaining terms that contain $-e_{1},-f_{1}$ or $-g_{1}$. Also,
since the image of $f''$ contains precisely one of $A_{1},B_{1},C_{1}$,
we conclude that the multiplicity of terms contain positive
$e_{1},f_{1},g_{1}$ must be zero. Thus, $\im(f'')$ is contained in the
configuration shown in \fullref{fig: no h's or e's}.

Now, note that in $\vd $ the sum of the multiplicities of the $e_{i}$'s is 
-2. This is true of the curve $A_{1}$ as well. Therefore the total multiplicity
of all other $e$ terms must vanish. But all other $e$ terms are of the form
$e_{i}-e_{i+1}$ or $e_{j}$. Since the former contribute nothing to the total
multiplicity, we conclude that there are no $e_{j}$ terms in the image of 
$f''$. Therefore $\im(f'')$ must be contained in the configuration shown in
\fullref{fig: disconnected X}.

\begin{figure}
\centerline{\includegraphics{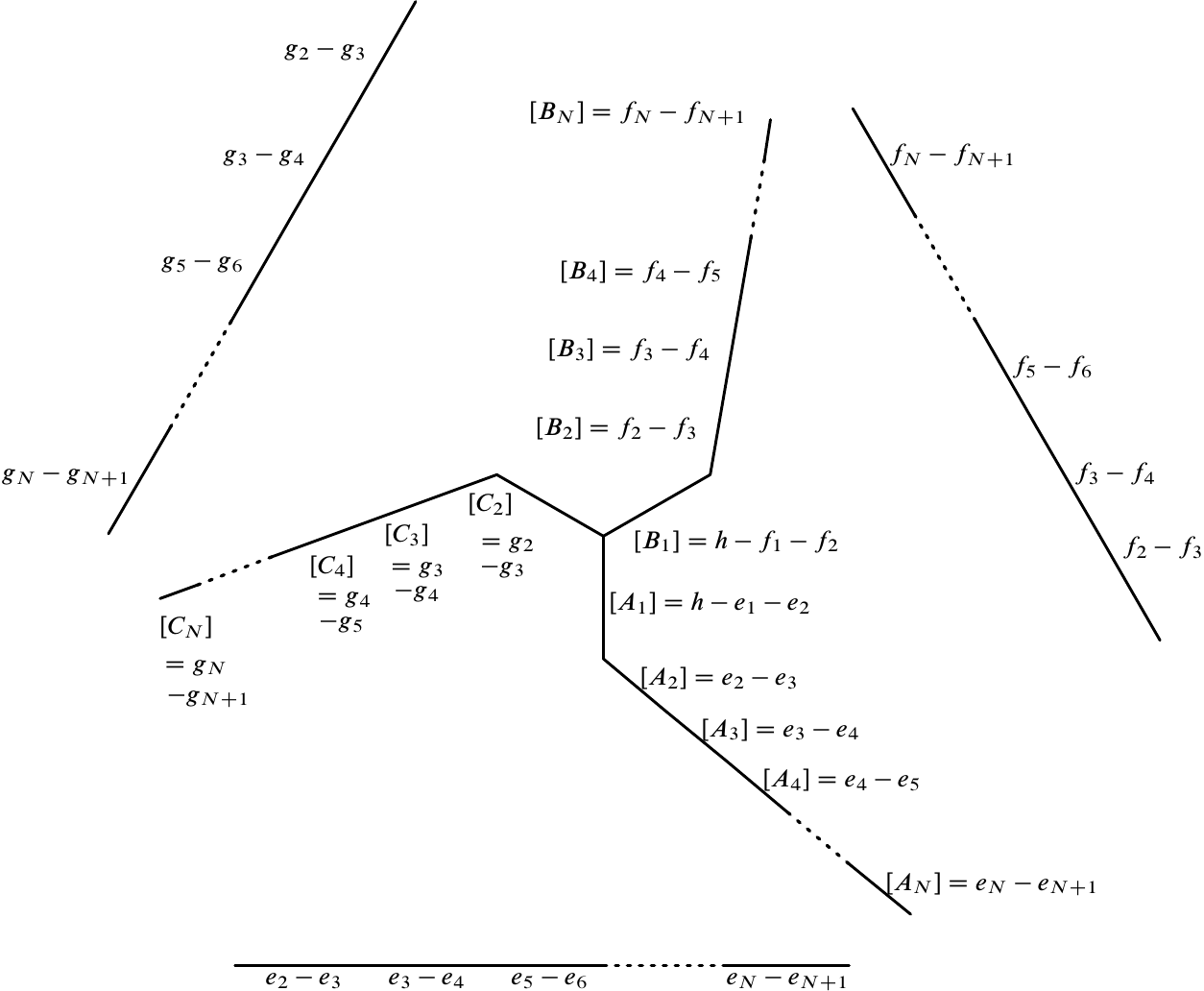}}
\caption{The remaining possible curves in $\im(f'')$}
\label{fig: disconnected X}
\end{figure}

But $\im(f'')$ is connected, and contains $h$ terms. Therefore it can not
contain nor be contained in any of the three outer parts of 
\fullref{fig: disconnected X}.
Therefore $\im(f'') \subset Y^{N}$. This contradicts 
our assumption, and therefore at least one of
$A_{1},B_{1},C_{1}$ is not in $\im(f'')$.

Without loss of generality, suppose $A_{1} \not\subset \im(f'')$. 
Let $d_{e,f},d_{e,g},d_{f,g}$ denote the
degree of $f''$ on the components $h-e_{1}-f_{1}, h-e_{1}-g_{1}, h-f_{1}-g_{1}$
respectively. Since $A_{1}$ is not contained in the image of $f''$, we must
have
\[
0< d_{e,f}+d_{e,g} \leq 3
\]
as these are the only multiplicities of $-e_{1}$ terms, and there are no
terms containing $-h$.

Furthermore, in order for $\im(f'')$ to simultaneously be connected and contain
$-e_{i}$ terms for $i>1$, it must be the case that $\im(f'')$ contains 
\emph{two} of 
\[
\{e_{1}, e_{1}-e_{2},e_{1}-\cdots -e_{N+1} \}.
\]
Thus
\[
d_{e,f}+d_{e,g} =3, \quad
d_{f,g} =0
\]
and $B_{1},C_{1} \not\subset \im(f'')$. This forces $\im(f'')$ to be contained
in the configuration shown in \fullref{fig: no dfg,b1,c1}.

\begin{figure}[ht!]
\centerline{\includegraphics{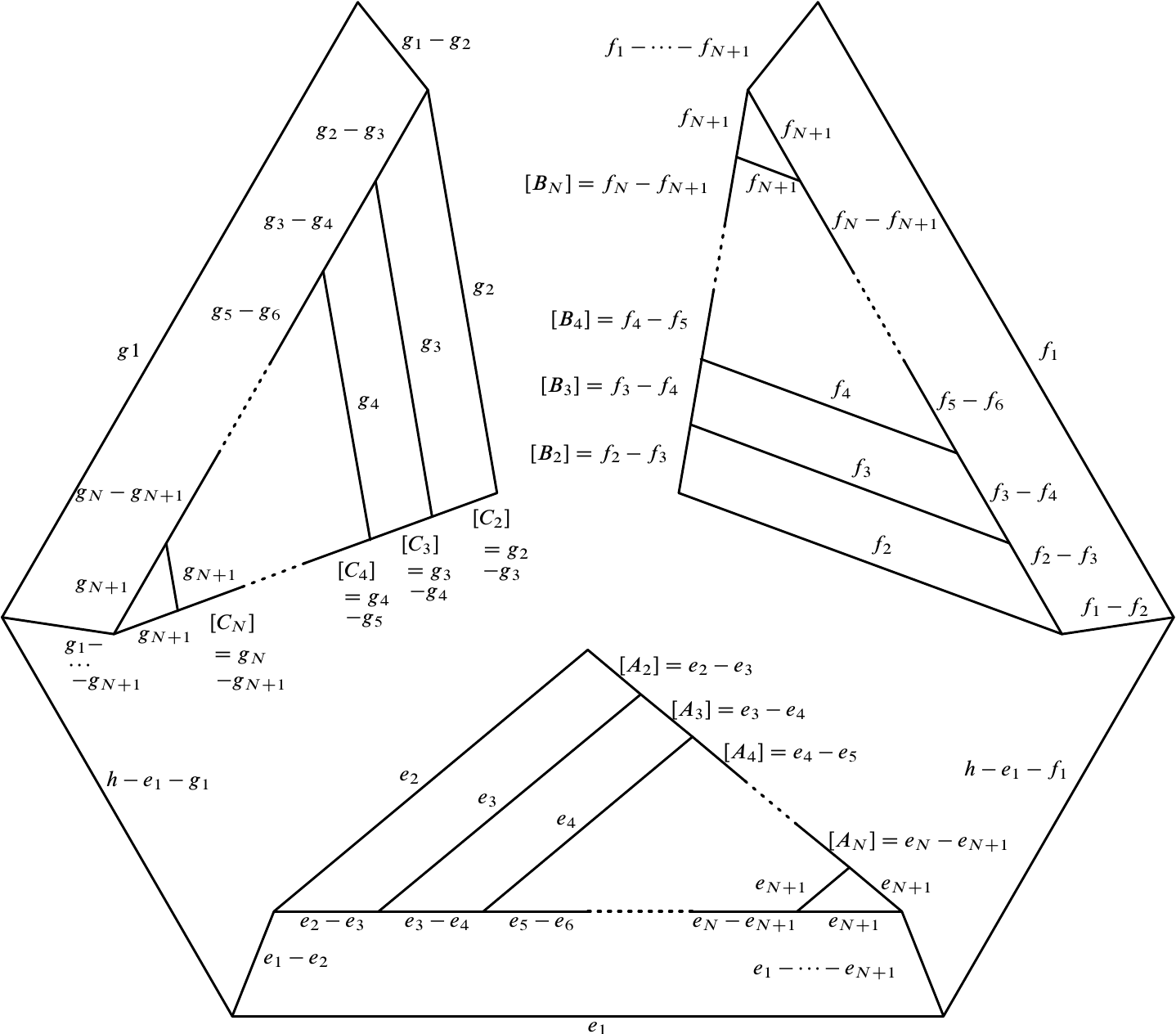}}
\caption{The other possibility for curves in $\im(f'')$}
\label{fig: no dfg,b1,c1}
\end{figure}

Again we have that $\im(f'')$ is connected and contains $-f_{i},-g_{j}$ for some
$i,j>1$. Therefore $\im(f'')$ contains at least one of 
$f_{1}-f_{2}, f_{1}-\cdots -f_{N+1}$ and also at least one of 
$g_{1}-g_{2}, g_{1}-\cdots -g_{N+1}$. But the multiplicity of $f_{1}$ and
$g_{1}$ in $\vd$ is $-1$. Therefore 
\[
d_{e,f},d_{e,g} \geq  2.
\]

This contradictions shows that our assumption $A_{1} \not\subset \im(f'')$
is incorrect. Therefore $A_{1} \subset \im(f'')$. An identical argument also
shows that $B_{1},C_{1} \subset\im(f'')$. However we showed above that
$A_{1},B_{1},C_{1} \not\subset \im(f'')$.

This contradiction shows that our original assumption is incorrect. Therefore
there does not exist a point $p \in \im(f'')$ such that $p \not\in Y^N$.
Thus $\im(f'') \subset Y^N$, and the result holds.
\end{proof}

\begin{rema}\label{remark: counter examples}
Note that this argument does not hold for general $a_{1},b_{1},c_{1}$.
For instance, it
is a fun exercise to show that there is more than one $\bT$--invariant
configuration of curves in $X$ in the following classes.
\begin{align*}
\beta_{1} =&2 (h-e_{1}-e_{2}) + (e_{2}-e_{3})\\
           & +2 (h-f_{1}-f_{2})+ (f_{2}-f_{3})\\
           & + 2 (h-g_{1}-g_{2}) + (g_{2}-g_{3})\\
\beta_{2}  =& 4 (h-e_{1}-e_{2}) + (e_{2}-e_{3}) + 2 (h-f_{1}-f_{2}) + 
               2 (h-g_{1}-g_{2}) + (g_{2}-g_{3})\\
\beta_{3} =& 4 (h-e_{1}-e_{2})+4 (e_{2}-e_{3})\\
           &+4 (h-f_{1}-f_{2})+4(f_{2}-f_{3})\\
           &+4 (h-g_{1}-g_{2}) +4 (g_{2}-g_{3}).
\end{align*}
\end{rema}

\subsection{A chain of rational curves}\label{section: local to global chain}
\begin{prop}\label{thm:chainIm}
Let $f\co  C\to X^{N+1}$ represent a point in 
$\Mbar_g(X^{N+1},\vd)$, where
$$
d_{1,1}>0,\quad  d_{2,j}=d_{3,j}=0, \quad j=1,\ldots,N.
$$
Then the image of $f$ is contained in 
the chain of rational curves 
\[
Y_A^N=A_1\cup\cdots\cup A_N
\]
defined in \fullref{sec:construction}.
\end{prop} 

Since $Y^N_A$ does not contain any of the curves $B_{i},C_{i}$, the blowups with
centers $p_{i}$ and $q_{i}$ in the construction of $X^{N+1}$ are extraneous. 
In order
to simplify the argument in this case, consider the space
\[
\dotsb \xrightarrow{\pi_{N+2}} \tX^{N+1} \xrightarrow{\pi_{N+1}} \tX^{N} 
\xrightarrow{\pi_{N}} \cdots  \xrightarrow{\pi_{1}} \bP^{3},
\]
where the construction of $\tX^{N+1}$ follows that of $X^{N+1}$, without the 
extraneous 
blowups. So $\tX^{i+1} \rightarrow \tX^{i}$ is the blowup of $\tX^{i}$ along
the point $p_{i}$, where $p_{i}$ is defined in 
\fullref{sec:construction}.  Thus, $\tX^{N+1}$ is deformation 
equivalent to the blowup of $\bP^{3}$ at $N+1$ points.
Since $Y^N_A$ does not contain the curves $B_{i},C_{i}$, clearly the formal
neighborhood of $Y^N_A$ in $\tX^{N+1}$ agrees with the construction in $X^{N+1}$.

\begin{figure}[ht!]
\centerline{\includegraphics{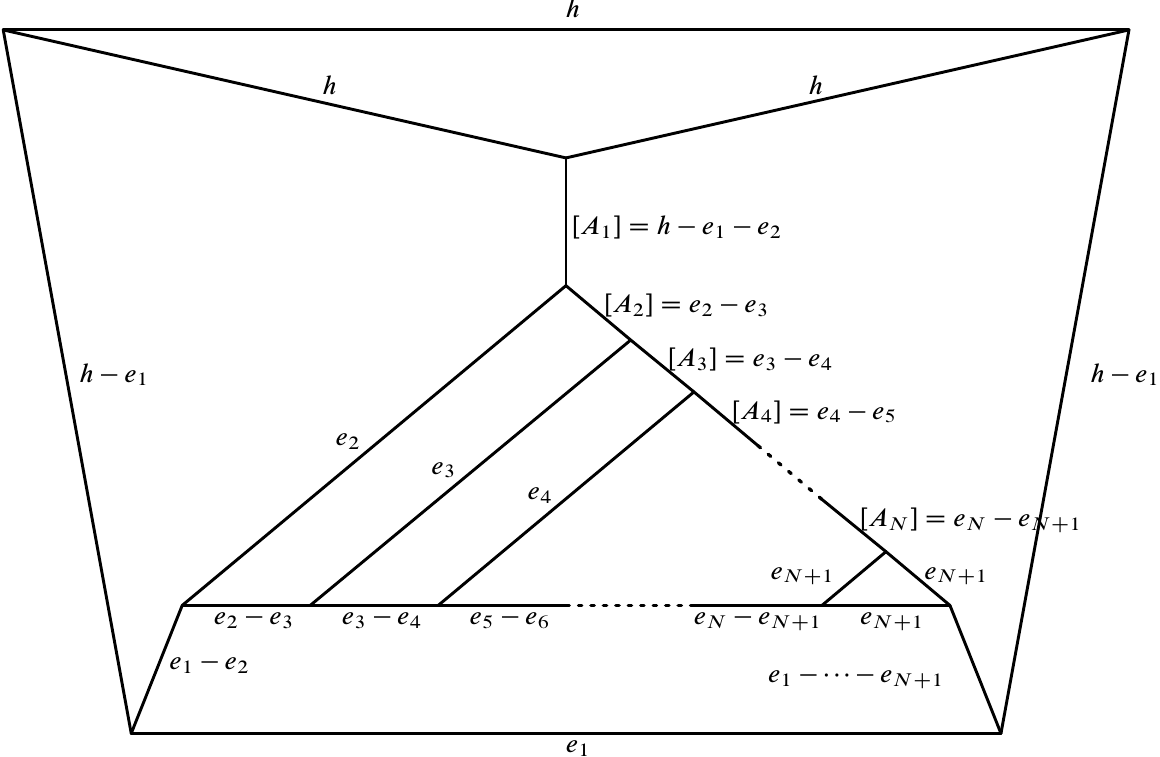}}
\caption{The $\bT$--invariant curves in $\tX^{N+1}$}
\label{fig: tX}
\end{figure}

We continue to let $E_{i}$ be the total transform of the exceptional divisor
over $p_{i}$, and $e_{i}$ be the class of a line in $E_{i}$. Furthermore,
we continue to let $H$ denote the pullback of the class of a hyperplane in
$\bP^{3}$, and $h$ be the class of a line in $H$. Then, $\{H,E_{i} \}$ is a
basis for $H_{4} (\tX^{N+1})$ and $\{h,e_{i} \}$ is a basis for 
$H_{2} (\tX^{N+1})$. The
non-zero intersection pairings are given as follows.
\begin{equation*}
\boxed{\begin{aligned}
	H \cdot H &= h  &H\cdot h &=pt\\ 
  E_{i}\cdot E_{i} &= -e_{i} &E_{i}\cdot e_{i} &=-pt\\
\end{aligned}}
\end{equation*}
The $\bT$--invariant curves in $\tX^{N+1}$ are shown together with their 
homology classes in \fullref{fig: tX}.

\begin{proof}[Proof of \fullref{thm:chainIm}]
As shown in above, we may use the toric
nature of $\tX^{N+1}$ to construct a stable map 
$[f''\co  \Sigma \rightarrow \tX^{N+1}] \in \Mbar_{g} (\tX^{N+1},\vd)$ such 
that 
$\im(f'' )$ is $\bT$--invariant, but $\im(f'') \not\subset Y^N_A$. We show that
this leads to a contradiction.

We study the class $f''_{*}[\Sigma] =\vd$. Note that the multiplicity of the
$-e_{1}$ term is the same as that of $h$. Furthermore, each $-e_{1}$ occurs
along with $h$, and there are no $-h$ terms. Therefore $\im(f'')$ can not
contain any terms containing positive $e_{1}$, nor can it contain any of
the curves in class $h$. Thus, the image of $f''$ is
contained in the configuration of curves shown in 
\fullref{fig: tX disconnected}.

\begin{figure}[ht!]
\centerline{\includegraphics{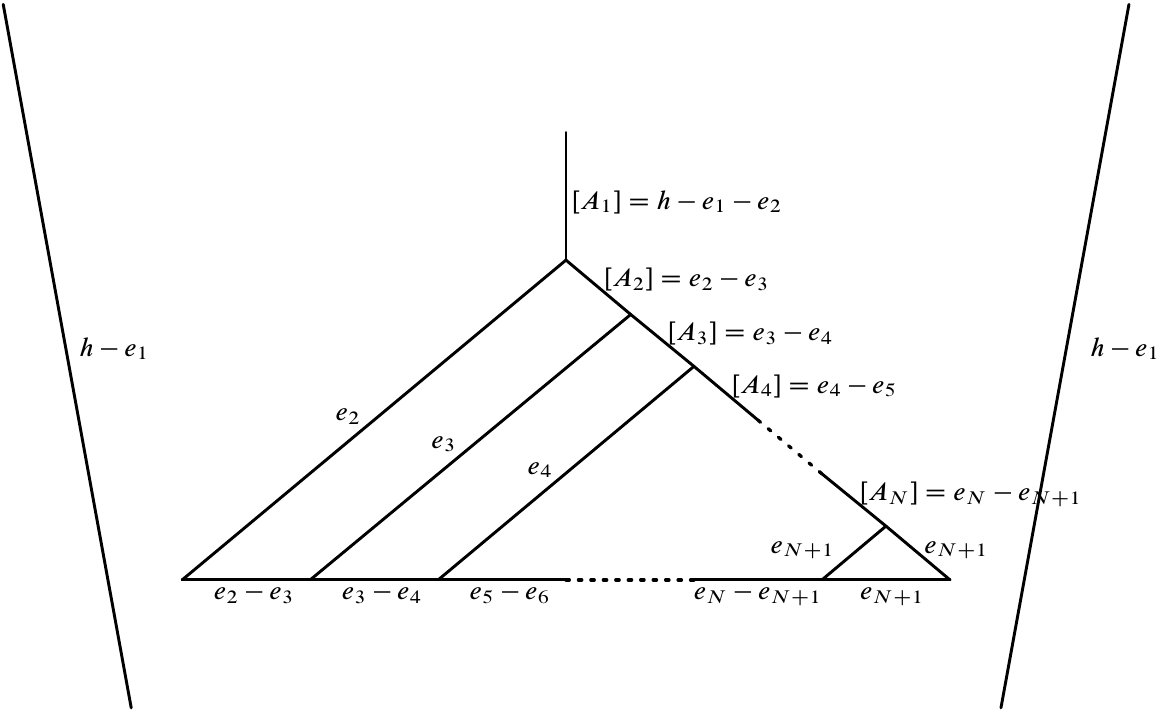}}
\caption{The possible curves in $\im(f'')$}
\label{fig: tX disconnected}
\end{figure}

Since $a_{1}>0$, it must be that $f''_{*}[\Sigma]$
contains at least one $e_{i}$ term with non-zero multiplicity for 
$i>1$. Also, $\im(f'')$ is connected and so
we conclude that the image of $f$ must not contain either 
of the curves of class $h-e_{1}$ in \fullref{fig: tX disconnected}.

Now, note that the total
multiplicity of the $e$ terms is $-2a_{1}$, and that the curve $A_{1}$ must 
also have this property. Therefore the sum of all other $e$ terms must be zero.
Since the other $e$ terms are of the form $e_{i}-e_{i+1}$ or $e_{j}$, we 
conclude that $\im(f'')$ does not contain any of the curves $e_{j}$. Thus
$\im(f'')$ is contained in the configuration depicted in 
\fullref{fig: tX complete}.

\begin{figure}[ht!]
\centerline{\includegraphics{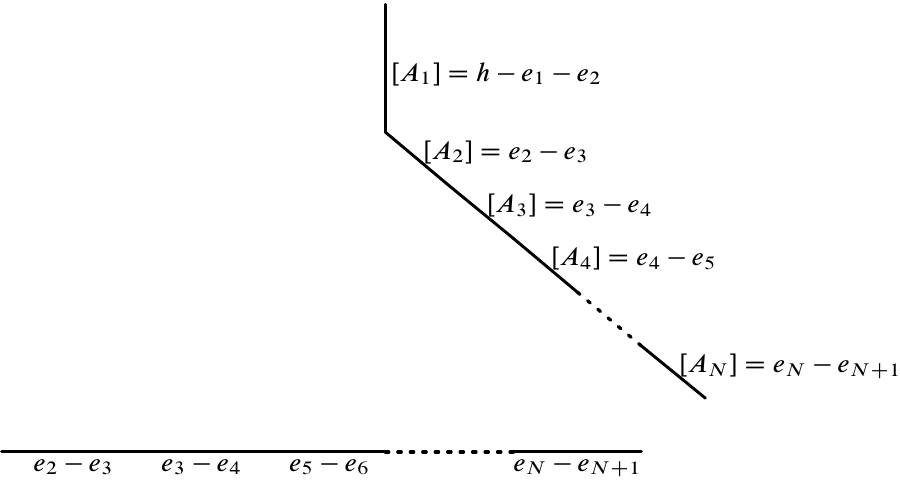}}
\caption{The remaining possible curves in $\im(f'')$}
\label{fig: tX complete}
\end{figure}

However, since $\tX^{N+1}$ is connected and contains $h$, we conclude that
$\im(f'') \subset Y^N_A$. This contradiction shows 
that our original assumption is incorrect, and that the result holds. 
\end{proof}

\section{Mathematical theory of the topological vertex}\label{sec:formal}
Let 
\begin{equation}\label{eqn:blowups}
 X^{N+1} \xrightarrow{\pi_{N+1}} X^{N} \xrightarrow{\pi_{N}}
\cdots \xrightarrow{\pi_{2}} X^{1} \xrightarrow{\pi_{1}} X^{0} =\bP^{3}
\end{equation}
be the toric blowups constructed in \fullref{sec:construction}.
Let $Y^N\subset X^{N+1}$ be the  minimal trivalent configuration, and let
$\hY^N$ be the formal completion of $X^{N+1}$ along $Y^N$.
Then $\hY^N$ is a nonsingular formal scheme, and
$\Mbar_g(\hY^N,\vd)$ is a separated formal
Deligne--Mumford stack with a perfect obstruction theory
of virtual dimension zero. It has a virtual fundamental
class when it is proper, which is not true in general.

\subsection{Formal Gromov--Witten invariants of $\hY^N$} \label{sec:cases}
In  \eqref{eqn:blowups}
$\bT=(\bC^\times)^3$ acts on $X^j$ and  the projections
$\pi_j$ are $\bT$--equivariant, so  $\hY^N$ is a formal scheme 
together with a $\bT$--action.  The point $s_0= A_1\cap B_1 \cap C_1$
is fixed by the $\bT$--action, so  $\bT$ acts on $T_{s_0} X^0$ and 
$\Lambda^3 T_{s_0} X^0$. Let $\bS$ be the rank 2 subtorus of $\bT$ which 
acts trivially on $\Lambda^3 T_{s_0} X^0$. The union of one dimensional
orbit closures of the $\bT$--action on $X^j$ is a configuration of
rational curves, which corresponds to a graph (see Figure 11).

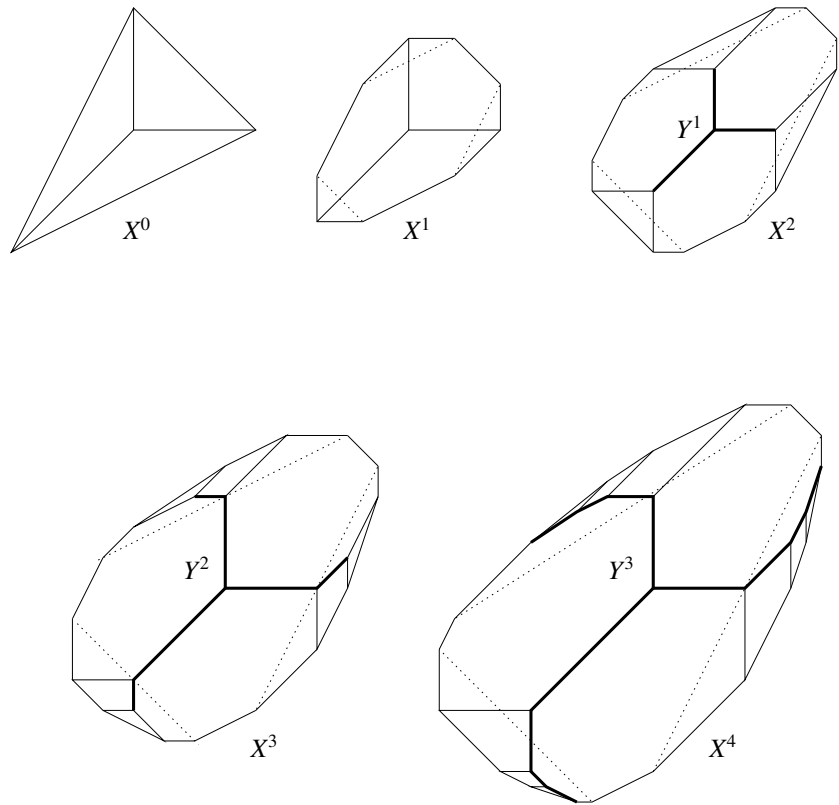
\begin{figure}[h]
\label{fig:XY}
\begin{center}
\labellist\small
\pinlabel {$X^0$} at 75 340
\pinlabel {$X^1$} at 240 340
\pinlabel {$X^2$} at 455 340
\pinlabel {$X^3$} at 150 35
\pinlabel {$X^4$} at 420 35
\pinlabel {$Y^1$} at 400 400
\pinlabel {$Y^2$} at 110 140
\pinlabel {$Y^3$} at 360 140
\endlabellist
\includegraphics[scale=0.65]{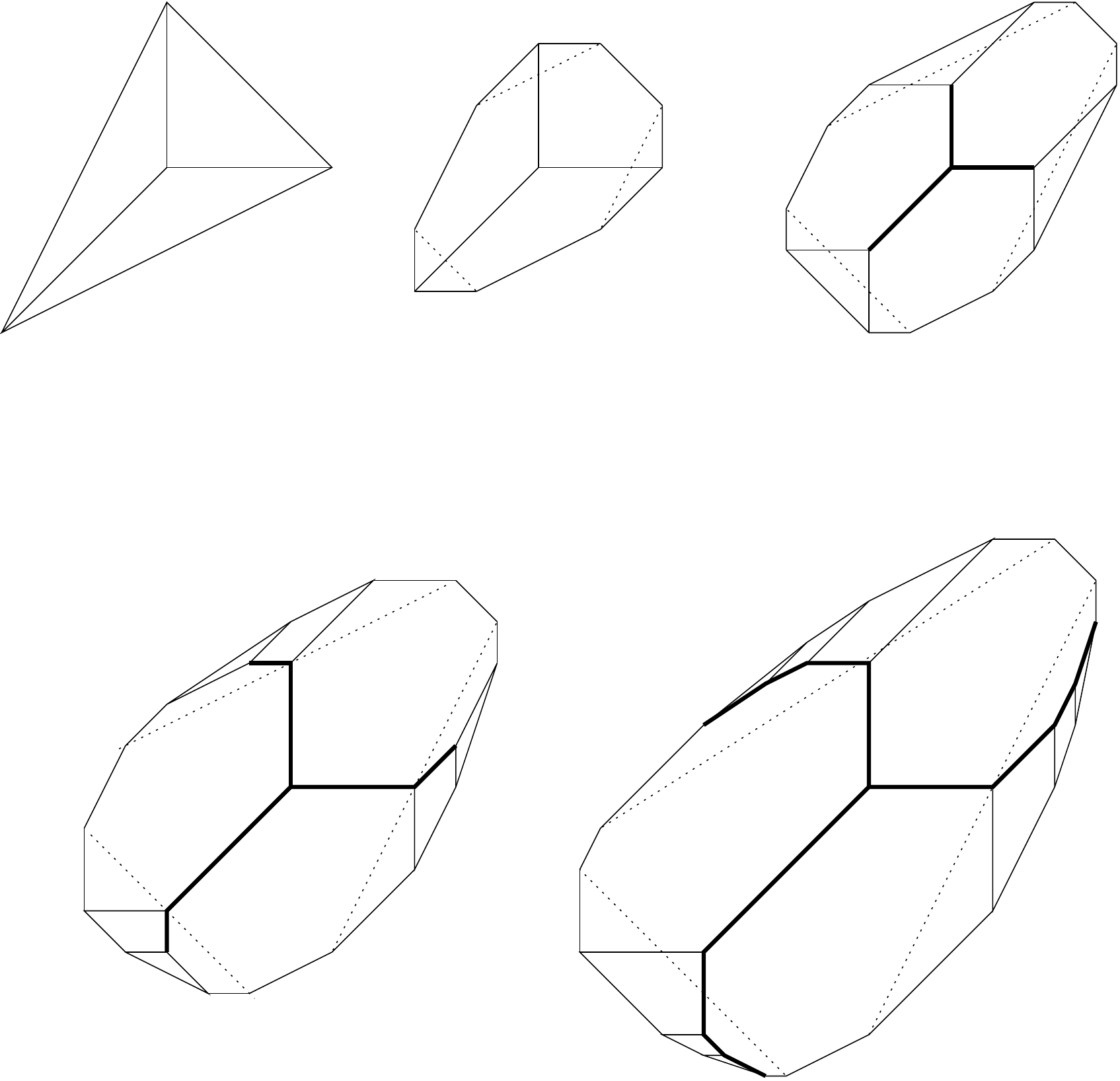}
\end{center}
\caption{Configurations of $\bT$--invariant curves}
\end{figure}

The $\bS$--action on $X^j$ can be read off from the slopes of the edges of the graph 
associated to $X^j$.More precisely, let $\Lambda_{\bS}=\Hom(\bS,\bC^\times)$ be
the group of irreducible characters of $\bS$. If we fix an identification 
$\bS\cong (\bC^\times)^2 $ then an element in $\Lambda_{\bS}$ is of the form 
$s_1^p s_2^q$ where $(s_1,s_2)$ are coordinates on $(\bC^\times)^2$ and 
$p,q\in \bZ$. The line segment associated to $C\cong \bP^1$ is
tangent to $(p,q)\in \bZ\oplus \bZ$ if the irreducible characters of the $\bS$--actions
on $T_x C$ and $T_y C$ are $s_1^p s_2^q$ 
and $s_1^{-p} s_2^{-q}$ (see Figure 12). Similarly, the $\bS$--action on $\hY^N$ can be read 
off from \fullref{fig:4cases} in \fullref{sec:computations}.

\begin{figure}[h]\label{fig:pq}
\begin{center}
\labellist\small
\pinlabel {$x$} [tr] at 10 10
\pinlabel {$y$} [bl] at 64 64
\pinlabel {$(p,q)$} [tl] at 20 20
\pinlabel {$(-p,-q)$} [tl] at 50 50
\endlabellist
\includegraphics[scale=0.65]{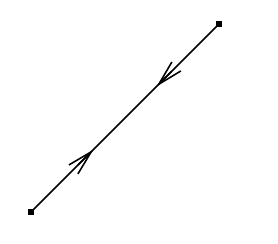}
\end{center}
\caption{The $\bS$--action can be read off from the slope}
\end{figure}
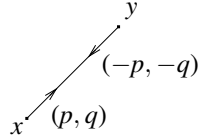

Let $u_1,u_2$ be a basis of $H^2_{\bS}(pt,\bZ)$ so that $H^2_{\bS}(pt,\bZ)=\bZ u_1\oplus \bZ u_2$.
For any nonzero effective class $\vd$ ($d_{i,j}\geq 0$), define
$$\tN^g_{\vd}(\hY^N)=\int_{[\Mbar_g(\hY^N,\vd)^{\bS}]^\vir}
  \frac{1}{e_{\bS}(N^\vir)}.$$
A priori $\tN^g_{\vd}(\hY^N)$ is a rational function in $u_1,u_2$
with $\bQ$ coefficients,  homogeneous of
degree 0. By results in Li--Liu--Liu--Zhou \cite{LLLZ},
$\tN^g_{\vd}(\hY^N)\in\bQ$ is a constant
function independent of $u_1,u_2$. We call $\tN^g_{\vd}(\hY^N)$ {\em formal 
Gromov--Witten invariants}. For the cases (i)--(iv) described in \fullref{sec:introduction},
\begin{eqnarray*}
N^g_{\vd}(Y^N)&=& \braket^{X^{N+1}}_{g,\vd}\\
&=&\int_{[\Mbar_g(X^{N+1},\vd)]^\vir} 1\\
&=&\int_{[\Mbar_g(X^{N+1},\vd)^{\bS} ]^\vir }
\frac{1}{e_{\bS}(N^\vir)}\\
&=&\int_{[\Mbar_g(\hY^N,\vd)^{\bS} ]^\vir }
\frac{1}{e_{\bS}(N^\vir)}\\
&=&\tN^g_{\vd}(\hY^N).
\end{eqnarray*}
As in \fullref{sec:introduction}, introduce formal variables $\la,t_{i,j}$, 
and define a generating function
\begin{equation}
F_N(\la;\vt)=\sum_{g\geq 0}\sum_\vd \lambda^{2g-2}\tN^g_\vd(\hY^N) e^{-\vd\cdot \vt}
\end{equation}
where
\[
\vt=(\lo{\vt}),\quad \vt_i = (t_{i,1},\ldots,t_{i,N}), \quad
\vd\cdot \vt=\sum_{i=1}^3\sum_{j=1}^N d_{i,j} t_{i,j}.
\]
The partition function of the formal Gromov--Witten invariants of $\hY^N$ is defined
to be
$$\tZ_N(\la;\vt)=\exp\left(F_N(\la;\vt)\right).$$
By connectedness and cyclic symmetry, we only need to
compute $\tN^g_\vd(\hY^N)$ in the following cases (see
\fullref{fig:4cases}):
\begin{enumerate}
\item[(D1)] $\vd=(\vd_1,\vo,\vo)$, 
      $d_{1,j}>0$ for $j\leq k$ and $d_{1,j}=0$ for $j>k$, 
      where $1\leq k\leq N$.
\item[(D2)] $\vd=(\vd_1,\vo,\vo)$, 
   $d_{1,j}>0$ for $k_1\leq j\leq k_2$ and $d_{1,j}=0$ otherwise, 
      where $2\leq k_1\leq k_2\leq N$.
\item[(D3)]  $\vd=(\vd_1,\vd_2,\vo)$, 
      $d_{1,m}>0$ for $m\leq j$ and $d_{1,m}=0$ for $m>j$, 
      $d_{2,m}>0$ for $m\leq k$ and $d_{2,m}=0$ for $m>k$, 
      where $1\leq j,k \leq N$.
\item[(D4)] $\vd=(\vd_1,\vd_2,\vd_3)$, 
      $d_{i,j}>0$ for $j\leq k_i$ and $d_{i,j}=0$ for $j>k_i$,  
      where $1\leq k_1,k_2,k_3 \leq N$.
\end{enumerate}
Any other $\tN^g_\vd(\hY^N)$ is either manifestly zero (because
$\Mbar_g(\hY^N,\vd)$ is empty) or is equal to one of the above case.
\begin{figure}[ht!]
\label{fig:4cases}
\begin{center}
\labellist\tiny
\pinlabel {(1)} [l] at 0 460
\pinlabel {$0$} [l] at 108 435
\pinlabel {$0$} [br] at 89 397
\pinlabel {$0$} [t] at 270 416
\pinlabel {$d_{1,1}$} [t] at 125 416
\pinlabel {$d_{1,2}$} [t] at 164 416
\pinlabel {$d_{1,k}$} [t] at 235 416
\pinlabel {(2)} [l] at 340 460
\pinlabel {$0$} [l] at 413 433
\pinlabel {$0$} [tl] at 398 400
\pinlabel {$0$} [t] at 430 416
\pinlabel {$0$} [t] at 502 416
\pinlabel {$0$} [t] at 647 416
\pinlabel {$d_{1,k_1}$} [t] at 544 416
\pinlabel {$d_{1,k_2}$} [t] at 614 416
\pinlabel {(3)} [l] at 0 190
\pinlabel {$0$} [l] at 107 270
\pinlabel {$0$} [tl] at 90 128
\pinlabel {$0$} [t] at 233 146
\pinlabel {$d_{2,k}$} [r] at 110 236
\pinlabel {$d_{2,1}$} [r] at 110 163
\pinlabel {$d_{1,1}$} [t] at 127 146
\pinlabel {$d_{1,j}$} [t] at 199 146
\pinlabel {(4)} [l] at 340 190
\pinlabel {$0$} [l] at 520 290
\pinlabel {$0$} [br] at 392 35
\pinlabel {$0$} [t] at 647 164
\pinlabel {$d_{3,k_3}$} [br] at 433 70
\pinlabel {$d_{3,1}$} [br] at 505 143
\pinlabel {$d_{2,1}$} [l] at 520 182
\pinlabel {$d_{2,k_2}$} [l] at 520 254
\pinlabel {$d_{1,1}$} [t] at 539 164
\pinlabel {$d_{1,k_1}$} [t] at 614 164
\endlabellist
\includegraphics[scale=0.4]{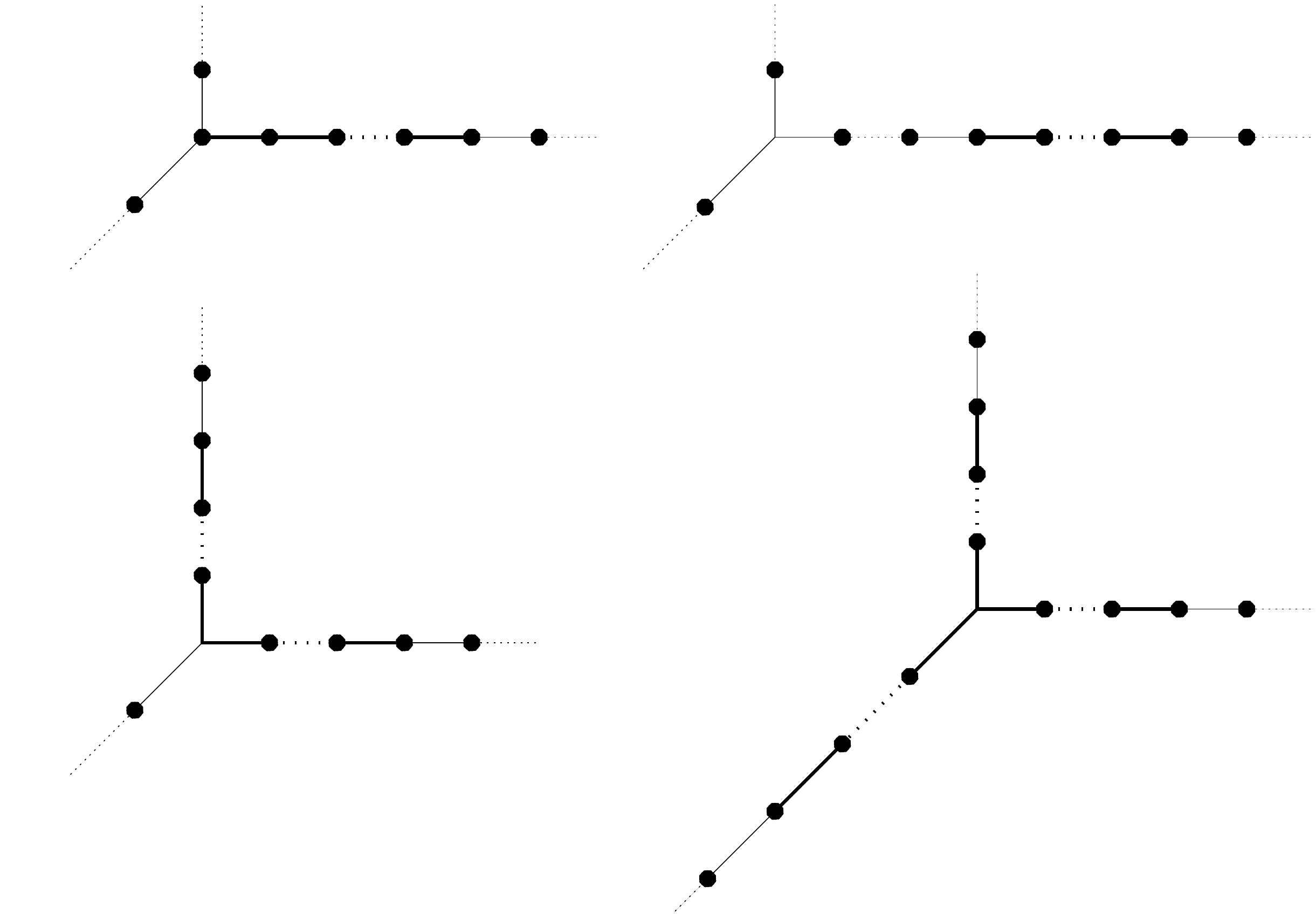}
\end{center}
\caption{Four cases}
\end{figure}
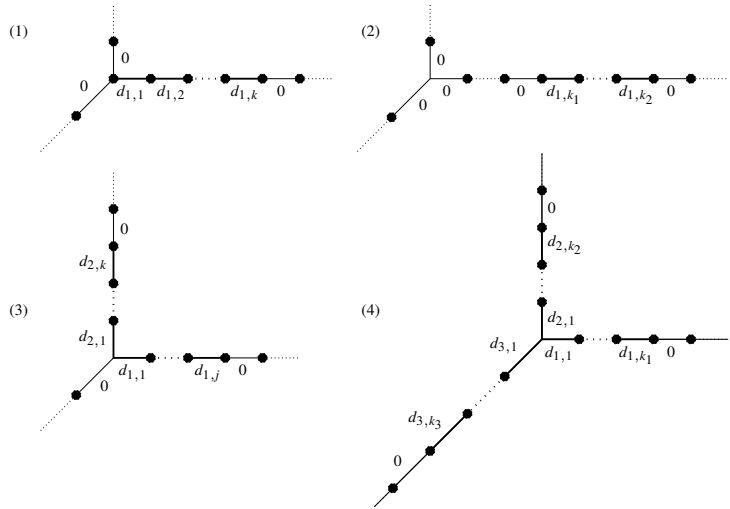
Let 
\[
F^1_N(\la;\vt_1),\quad F^2_N(\la;\vt_1),\quad
F^3_N(\la;\vt_1,\vt_2),\quad F^4_N(\la;\vt)\quad
\]
denote the contribution to $F_N(\la;\vt)$ from (D1), (D2), (D3), (D4),
respectively.
Then
\begin{eqnarray*}
F_N(\la;\vt)&=&\sum_{i=1}^3 F^1_N(\la;\vt_i) +\sum_{i=1}^3 F^2_N(\la;\vt_i) \\
&&+F^3_N(\la;\vt_1,\vt_2)+ F^3_N(\la;\vt_2,\vt_3)+  F^3_N(\la;\vt_3,\vt_1)
+F^4_N(\la;\vt).
\end{eqnarray*}

\subsection{Summary of results}
Let $C_g$ be defined by
$$
\sum_{g=0}^\infty C_g t^{2g}=\left(\frac{t/2}{\sin(t/2)}\right)^2
$$
as before. Then
$$
\sum_{g\geq 0} C_g n^{2g-3} \lambda^{2g-2}=\frac{-1}{n[n]^2}
$$
where 
$$
[n]=q^{n/2}-q^{-n/2},\quad q=e^{\sqrt{-1}\la}.
$$

In \fullref{sec:computations}, we will compute $\tN^g_\vd(\hY^N)$ by
the mathematical theory of the topological vertex and obtain the following
results. We will do the computations by the physical theory
of the topological vertex in \fullref{sec:vertex}.

\begin{prop}\label{thm:Fone}
Suppose that
$$
\vd_1=(d_1,\ldots,d_N),\quad \vd_2=\vd_3=(0,\ldots,0).
$$
where $d_1>0$. Then
\[
\tN^g_{\vd}(\hY^N)=\left\{ \begin{array}{ll}
C_g n^{2g-3} & 
\begin{array}{l}d_1 = d_2 = \cdots = d_k = n > 0\textit{ and}\\
 d_{k+1} = d_{k+2} = \cdots = d_N=0 \textit{ for some }1\leq k\leq N
\end{array}\\
 0  & \textit{otherwise}   
\end{array}\right.
\]
which is equivalent to
\begin{equation}\label{eqn:Fone}
F_N^1(\la;\vt_1)=\sum_{n > 0} \frac{-1}{n[n]^2}\sum_{1\leq k\leq N}
e^{-n(t_{1,1}+\cdots+ t_{1,k})}
\end{equation}
\end{prop}

\fullref{thm:Fone} is equivalent to \fullref{thm:chain}.

\begin{prop}\label{thm:Ftwo}
Suppose that
$$
\vd_1=(d_1,\ldots,d_N),\quad \vd_2=\vd_3=(0,\ldots,0).
$$
where $d_1=0$. Then
\[
\tN^g_{\vd}(\hY^N)=\left\{ \begin{array}{ll}
-C_g n^{2g-3} &  \begin{array}{ll}d_j = n > 0\textit{ for } k_1\leq j\leq k_2 \textit{ and}\\
d_j=0 \textit{ otherwise, where } 2\leq k_1\leq k_2\leq N\end{array} \\
 0  & \textit{otherwise}   
\end{array}\right.
\]
which is equivalent to
\begin{equation}\label{eqn:Ftwo}
F_N^2(\la;\vt_1)=\sum_{n > 0} \frac{1}{n[n]^2}\sum_{2\leq k_1\leq k_2\leq N}
e^{-n(t_{1,k_1}+\cdots+ t_{1,k_2})}
\end{equation}
\end{prop}

\fullref{thm:Ftwo} corresponds to a chain of $(0,-2)$ rational curves.

\begin{prop}\label{thm:Fthree}
Suppose that 
$$
\vd_1=(n,0,\ldots,0),\quad  d_{2,1}>0,\quad \vd_3=(0,\ldots,0)
$$
where $n>0$. Then
\[
\tN^g_{\vd}(\hY^N)=\left\{\begin{array}{ll}
-C_g n^{2g-3} & \begin{array}{l}d_{2,1} = d_{2,2}=\cdots = d_{2,k} = n \textit{ and}\\
 d_{k+1} = d_{k+2} = \cdots = d_N=0 \textit{ for some }1\leq k\leq N
\end{array}\\
0 &\textit{otherwise}.
\end{array}\right.
\]
\end{prop}

\begin{prop}\label{thm:decreasing}
Suppose that $d_{i,1}>0$.
Then $\tN^g_{\vd}(\hY^N)=0$ unless
$$
d_{i,1}\geq d_{i,2}\geq \cdots \geq d_{i,N}.
$$
\end{prop}

\begin{prop}\label{thm:Ffour}
Suppose that 
\[
d_{i,j}=\left\{\begin{array}{ll}
d_i>0 & j\leq k_i\\
0 &j>k_i\end{array}\right.
\]
where $i=1,2,3$ and $1\leq k_i\leq N$. Then
\begin{enumerate}
\item[(a)] 
$
\tN^g_{\vd_1,\vd_2,\vo}(\hY^N)=\left\{\begin{array}{ll}
-C_g n^{2g-3} & d_1=d_2=n>0\\
0 & \textit{otherwise} \end{array}\right.
$
\medskip
\item [(b)]
$
\tN^g_{\vd_1,\vd_2,\vd_3}(\hY^N)=\left\{\begin{array}{ll}
C_g n^{2g-3} & d_1=d_2=d_3=n>0\\
0 & \textit{otherwise} \end{array}\right.
$
\end{enumerate}
\end{prop}
\fullref{thm:decreasing} and \fullref{thm:Ffour}
are consistent with \fullref{thm:tri}.

Let $N=1$ in \fullref{thm:Fone} and \fullref{thm:Ffour},
we get
\begin{coro}\label{thm:tZone-Q}
$$
\tZ_1(\la;\vt)=\exp\left(\sum_{n=1}^\infty
\frac{Q_n(\vt) }{-n[n]^2}\right)
$$
where $\vt=(\lo{t})$ and 
$$
Q_n(\vt)=e^{-nt_1}+e^{-nt_2}+ e^{-nt_3}-
e^{-n(t_1+t_2)}-e^{-n(t_2+t_3)}-e^{-n(t_3+t_1)}
+e^{-n(t_1+t_2+t_3)}.
$$
\end{coro}

Finally, we will derive the following expression
of $\tZ_N(\la;\vt)$, where the notation is the same as that 
 in \fullref{sec:introduction}:
\begin{prop}\label{thm:tZ}
\[
\tZ_N(\la;\vt)=\exp\biggl(\sum_{i=1}^3 F_N^2(\la;\vt_i)\biggr) 
\sum_{\vmu} \tilde{\cW}_\vmu(q) \prod_{i=1}^3
(-1)^{|\mu^i|} e^{-|\mu^i|t_{i,1}} s_{(\mu^i)^t}(u^i(q,\vt_i)).
\]
where 
\[
F_N^2(\lam;\vt_i)=\sum_{n=1}^\infty \frac{1}{n[n]^2}\sum_{2\leq k_1\leq k_2 \leq N}
e^{-n(t_{i,k_1}+\cdots +t_{i,k_2}) }
\]
\end{prop}
In particular, when $N=1$ we have
$$F_1^2(\lam;\vt_i)=0, \quad
  s_{(\mu^i)^t}(u^i(q,\vt_i))=\cW_{(\mu^i)^t}(q),$$
which gives the following corollary.
\begin{coro}\label{thm:tZone-W}
\[
\tZ_1(\la;\vt)=\sum_{\vmu} \tilde{\cW}_\vmu(q) \prod_{i=1}^3
(-1)^{|\mu^i|} e^{-|\mu^i|t_i} \cW_{(\mu^i)^t}(q).
\]
where $\vt=(\lo{t})$.
\end{coro}

Equation \eqref{eqn:Zone} in \fullref{sec:introduction} 
follows from \fullref{thm:tZone-Q} and \fullref{thm:tZone-W}.

\subsection{Three-partition Hodge integrals}
\label{sec:Hodge}
Three-partition Hodge integrals arise when we calculate 
$$
\tN^g_{\vd}(\hY^N)=\int_{[\Mbar_g(\hY^N,\vd)^{\bS}]^{\vir} } \frac{1}{e_{\bS}(N^\vir)}
$$
by virtual localization (see Li--Liu--Liu--Zhou \cite[Section 7]{LLLZ} for such calculations). 
We recall their definition in this subsection.

Let $\lo{w}$ be formal variables, where $w_3=-w_1-w_2$. Let $w_4=w_1$. Write $\bw=\three{w}$.
For $\vmu=(\up{\mu})\neq (\emptyset,\emptyset,\emptyset)=\vec{\emptyset}$, define
$$
d^1_\vmu=0,\quad
d^2_\vmu=\ell(\mu^1),\quad
d^3_\vmu=\ell(\mu^1)+\ell(\mu^2),\quad
\ell(\vmu)=\sum_{i=1}^3\ell(\mu^i).
$$
The {\em three-partition Hodge integrals} are defined by
\begin{multline}\label{eqn:Gg}
G_{g,\vmu}(\bw)=\frac{ (-\sqrt{-1})^{\ell(\vmu)} }{|\Aut(\vmu)|}
 \prod_{i=1}^3\prod_{j=1}^{\ell(\mu^i)}
\frac{\prod_{a=1}^{\mu^i_j-1}(\mu^i_j w_{i+1} + a w_i) }
     {(\mu^i_j-1)!w_i^{\mu^i_j-1} }\\
\cdot \int_{\Mbar_{g,\ell(\vmu)} }
\prod_{i=1}^3\frac{\Lambda_g^\vee(w_i)w_i^{\ell(\vmu)-1} }
{\prod_{j=1}^{\ell(\mu^i)}(w_i(w_i-\mu^i_j\psi_{d^i_\vmu+j}))}
\end{multline}
where
$$
\Lambda_g^\vee(u)=u^g-\lam_1 u^{g-1}+\cdots + (-1)^g\lam_g.
$$
Note that $G_{g,\vmu}(\lo{w})$ has a pole along $w_i=0$ if $\mu^i\neq
\emptyset$.
The following cyclic symmetry is clear from the definition:
\begin{equation}
\begin{aligned}
 G_{g,\up{\mu}}(\lo{w})
&=G_{g,\mu^2,\mu^3,\mu^1}(w_2,w_3,w_1) \\
&=G_{g,\mu^3,\mu^1,\mu^2}(w_3,w_1,w_2)
\end{aligned}
\end{equation}

Note that
$$
\sqrt{-1}^{\ell(\vmu)}G_{g,\vmu}(\bw)\in \QQ\three{w}
$$
is homogeneous of degree $0$, so
$$
G_{g,\vmu}(w_1,w_2,-w_1-w_2)
=G_{g,\vmu}(1,\frac{w_2}{w_1}, -1-\frac{w_2}{w_1}).
$$
Introduce variables $\lam$, $p^i=(p^i_1,p^i_2,\ldots)$, $i=1,2,3$.
Given a partition $\mu$,
define
$$
p^i_\mu=p^i_1\cdots p^i_{\ell(\mu)}
$$
for $i=1,2,3$. In particular, $p^i_\emptyset=1$.
Write
$$
\bp=(\up{p}),\ \ \
\bp_\vmu=p^1_{\mu^1}p^2_{\mu^2}p^3_{\mu^3}.
$$
Define generating functions
\begin{eqnarray*}
G_\vmu(\lam;\bw)&=&\sum_{g=0}^\infty \lam^{2g-2+\ell(\vmu)}G_{g,\vmu}(\bw),\\
G(\lam;\bp;\bw)&=&\sum_{\vmu\neq \vec{\emptyset}} G_\vmu(\lam;\bw)\bp_\vmu,\\
G^\bu(\lam;\bp;\bw)&=&\exp(G(\lam;\bp;\bw))
=1+\sum_{\vmu\neq \vec{\emptyset} } G^\bu_\vmu(\lam;\bw)\bp_\vmu.
\end{eqnarray*}
 
In particular, 
\[
G_{g,\mu,\emptyset,\emptyset}(1,0,-1)=\left\{\begin{array}{ll}
-\sqrt{-1}d^{2g-2}b_g, & \mu=(d),\\
0, & \ell(\mu)>1.\end{array}\right.
\]
where 
$$
b_g=\left\{\begin{array}{ll}
1, &g=0,\\
\int_{\Mbar_{g,1}} \psi_1^{2g-2}\lambda_g & g>0. 
\end{array}\right.
$$
It was proved by Faber and Pandharipande \cite{FP} that
\begin{equation}
\sum_{g=0}^\infty b_g t^{2g}=\frac{t/2}{\sin(t/2)}.
\end{equation}
So
\begin{equation}\label{eqn:Gzero}
G_{(n),\emptyset,\emptyset}(\lambda;1,0,-1)=\frac{-\sqrt{-1}}{2n\sin(\la n/2)}
=\frac{1}{n[n]}.
\end{equation}
Similarly, we have
\begin{equation}\label{eqn:Gone}
G_{(n),\emptyset,\emptyset}(\lambda;1,-1,0)=\frac{(-1)^{n-1}}{n[n]}.
\end{equation}
The following formula of three-partition Hodge integrals
was derived in Li--Liu--Liu--Zhou \cite{LLLZ}:
\begin{equation}\label{eqn:Gthree}
G^\bu_\vmu(\la;\bw)=\sum_{|\nu^i|=|\mu^i|}\prod_{i=1}^3
\left(q^{\frac{1}{2}\kappa_{\nu^i}\frac{w_{i+1}}{w_i} } 
\frac{\chi_{\nu^i}(\mu^i) }{z_{\mu^i}} \right)\tilde{\cW}_{\vnu}(q),
\end{equation}

In particular, 
\begin{eqnarray}
\label{eqn:MarinoVafa}
G^\bu_{\mu,\emptyset,\emptyset}(\la;\!1,\!\tau,\!-\tau{-}1)
&=&\!\!\sum_{|\nu|=|\mu|}\!\! \frac{\chi_\nu(\mu)}{z_\mu}
  q^{\frac{1}{2}\kappa_\nu \tau}\cW_\nu(q).\\
\label{eqn:Gtwo}
\qquad G^\bu_{\mu^1,\mu^2,\emptyset}(\la;\!1,\!\tau,\!-\tau{-}1)
&=&\!\!\!\!\sum_{|\nu^i|=|\mu^i|}\!\!\!\!
  \frac{\chi_{\nu^1}(\mu^1)}{z_{\mu^1}}
  \frac{\chi_{\nu^2}(\mu^2)}{z_{\mu^2}}
  q^{\frac{1}{2}\left(\kappa_{\nu^1}\tau -\kappa_{\nu^2}\tau^{-1}\right)}
  \cW_{\nu^1 (\nu^2)^t }(q).
\end{eqnarray}
Equation \eqref{eqn:MarinoVafa} is equivalent to the formula of one-partition
Hodge integrals conjectured in Mari\~no--Vafa \cite{MV}, which was proved
in Liu--Liu--Zhou \cite{LLZ03}
and Okounkov--Pandharipande \cite{OP}. Equation \eqref{eqn:Gtwo} is equivalent to the formula
of two-partition Hodge integrals proved in Liu--Liu--Zhou \cite{LLZ}.

\subsection{Relative formal GW invariants of the topological vertex}
Symplectic relative Gromov--Witten theory was developed by Li and Ruan
\cite{LR}, and Ionel and Parker \cite{IP03,IP04}.
The mathematical theory of the topological vertex in \cite{LLLZ} is  
based on Jun Li's algebraic relative Gromov--Witten theory \cite{Li01,Li02}.

Given a triple of partitions
$\vmu=(\up{\mu})\neq \vec{\emptyset}$
and a triple of integers $\vn=(\lo{n})$, let 
$$
F^\bu_{\chi,\vmu}(\vn)\in \QQ
$$
be the disconnected relative formal GW invariants of the 
topological vertex defined in \cite{LLLZ}. Introduce variables 
$\la, p^i_j$ as in \fullref{sec:Hodge},
and define generating functions
\begin{eqnarray*}
F_{\vmu}^\bu(\la;\vn)&=&
\sum_{\chi} \la^{-\chi+\ell(\vmu)}F^\bu_{\chi,\vmu}(\vn)\\
F^\bu(\la;\vp;\vn)&=& 
1+\sum_{\vmu\neq \vec{\emptyset}} F^\bu_{\vmu}(\la;\vn) \bp_{\vmu}\\
F(\la;\vp;\vn)&=&\log(F^\bu(\la;\vp;\vn))
=\sum_{\vmu\neq \vec{\emptyset}} F_{\vmu}(\la;\vn)\vp_{\vmu}\\
F_{\vmu}(\la;\vn)
&=& \sum_{g=0}^\infty \la^{2g-2+\ell(\vmu)}F_{g,\vmu}(\vn).
\end{eqnarray*}
By virtual localization, $F_{g,\vmu}(\vn)$ can be expressed in terms of
three-partition Hodge integrals and double Hurwitz numbers. We have
\begin{multline}\label{eqn:FG}
F^\bu_\vmu(\lam;\vn)
=\\
(-1)^{\sum_{i=1}^3 (n_i-1)|\mu^i| }(-\sqrt{-1})^{\ell(\vmu)}
\!\!\!\!
\sum_{|\nu^i|=|\mu^i|}
\!\!\!\!
G^\bu_{\vnu}(\la;\bw)
\prod_{i=1}^3 z_{\nu^i}
\Phi^\bu_{\nu^i,\mu^i}\Bigl(\sqrt{-1}\bigl(n_i{-}\frac{w_{i+1}}{w_i}\bigr)\Bigr)
\end{multline}
where 
\begin{equation} \label{eqn:burnside}
\Phi^\bu_{\nu,\mu}(\lam)=\sum_\chi H^\bu_{\chi,\nu,\mu}
\frac{\la^{-\chi+\ell(\nu)+\ell(\mu)} }{(\chi+\ell(\nu)+\ell(\mu))!}
=\sum_{\eta}e^{\kappa_\eta \lam/2}
\frac{\chi_\eta(\nu) }{z_{\nu} }
\frac{\chi_\eta(\mu) }{z_{\mu} }.
\end{equation} 
is a generation function of disconnected double 
Hurwitz numbers $H^\bu_{\chi,\mu,\nu}$.
Equations \eqref{eqn:Gthree}, \eqref{eqn:FG}, and \eqref{eqn:burnside}
imply 
\begin{equation}\label{eqn:FW}
F^\bu_\vmu(\la;\vn)=\sum_{|\nu^i|=|\mu^i|}\prod_{i=1}^3
\left(q^{\frac{1}{2}\kappa_{\nu^i} n_i} 
\frac{\chi_{\nu^i}(\mu^i) }{z_{\mu^i}} \right)\tilde{\cW}_{\vnu}(q),
\end{equation}
Note that the $\bw$ dependence on the right hand side of \eqref{eqn:FG} cancels and 
the right hand side of \eqref{eqn:FW} is independent of $\bw$. Since 
$\Phi^\bu_{\nu\mu}(0)=\delta_{\nu\mu}/z_\mu$, we have
$$
F^\bu_{\mu}(\frac{w_2}{w_1},\frac{w_3}{w_2},\frac{w_1}{w_3})=
(-1)^{\sum_{i=1}^3 (n_i-1)|\mu^i| }(-\sqrt{-1})^{\ell(\vmu)}
\sum_{|\nu^i|=|\mu^i|} G^\bu_{\vnu}(\la;\lo{w}).
$$
Also
$F^\bu_\vmu(\la;\vn)$ is independent of $n_i$ if $\mu^i$ is empty. 
\begin{eqnarray}
\label{eqn:gmu}
F_{g,\mu,\emptyset,\emptyset}(0,n_2,n_3)
  &=&(-1)^{|\mu|}(-\sqrt{-1})^{\ell(\mu)}
  G_{g,\mu,\emptyset,\emptyset}(1,0,-1) \\
\label{eqn:gmunu}
F_{g,\mu,\nu,\emptyset}(-1,0,n_3)
  &=&(-1)^{|\nu|}(-\sqrt{-1})^{\ell(\mu)+\ell(\nu)}
  G_{g,\mu,\nu,\emptyset}(1,-1,0)
\end{eqnarray}
From \eqref{eqn:Gzero}, \eqref{eqn:Gone}, \eqref{eqn:gmu}, \eqref{eqn:gmunu}
we conclude that if $\mu\neq \emptyset$ then
\begin{eqnarray}
F_{\mu,\emptyset,\emptyset}(\la;0,n_2,n_3)&=&
  \left\{\begin{array}{ll}
  \frac{(-1)^{n-1}\sqrt{-1}}{n[n]} & \mu=(n)\\
  0 & \ell(\mu)>1.
  \end{array}\right. \\
F_{\mu,\emptyset,\emptyset}(\la;-1,n_2,n_3)&=&
  \left\{\begin{array}{ll}
  \frac{(-1)^n\sqrt{-1}}{n[n]} & \mu=(n)\\
  0 & \ell(\mu)>1.
  \end{array}\right.
\end{eqnarray}
If $\mu\neq \emptyset, \nu\neq\emptyset$ then
(see Liu--Liu--Zhou \cite[p7]{LLZ} for details)
\begin{equation}
F_{\mu,\nu,\emptyset}(\la;-1,0,n_3)=
\left\{\begin{array}{ll}
\frac{1}{n} & \mu=\nu=(n)\\
0 & \textup{otherwise}.\end{array}\right.
\end{equation}
We also have
$$F_{(1),(1),\emptyset}(\la;0,0,0)=-1,\ \ \ 
  F_{(1),(1),(1)}(\la;0,0,0)=-\sqrt{-1}[1].$$

\subsection{Computations}\label{sec:computations}

The $\bS$--action on $\hY^N$ can be read off from \fullref{fig:4cases}
as explained in the first two paragraphs of \fullref{sec:cases}.

\begin{figure}[hb!]\label{fig:minimal}
\labellist\small
\hair=1.5pt
\pinlabel {$d_{1,1}$} [t] at 180 144
\pinlabel {$d_{1,2}$} [br] at 220 162
\pinlabel {$d_{1,3}$} [br] at 249 213
\pinlabel {$d_{1,4}$} [r] at 280 280
\pinlabel {$d_{2,1}$} [l] at 162 162
\pinlabel {$d_{2,2}$} [t] at 142 180
\pinlabel {$d_{2,3}$} [tl] at 94 165
\pinlabel {$d_{2,4}$} [tl] at 30 130
\pinlabel {$d_{3,1}$} [tl] at 142 126
\pinlabel {$d_{3,2}$} [l] at 126 92
\pinlabel {$d_{3,3}$} [bl] at 140 56
\pinlabel {$d_{3,4}$} [bl] at 192 22
\endlabellist
\begin{center}
\includegraphics[scale=0.6]{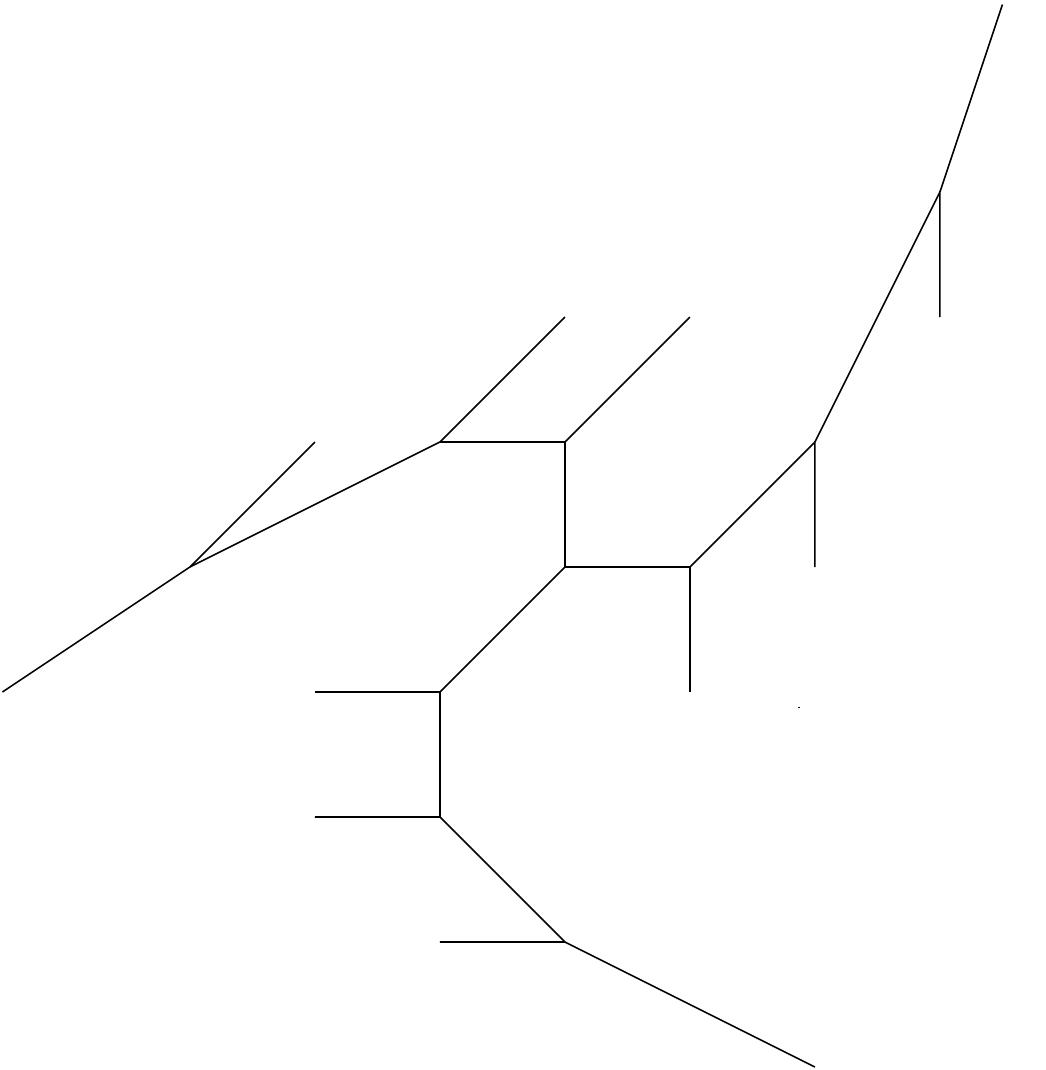}
\end{center}
\caption{Graph of $\hY_N$}
\end{figure}

\begin{figure}[hb!]\label{fig:break}
\begin{center}
\includegraphics[scale=0.6]{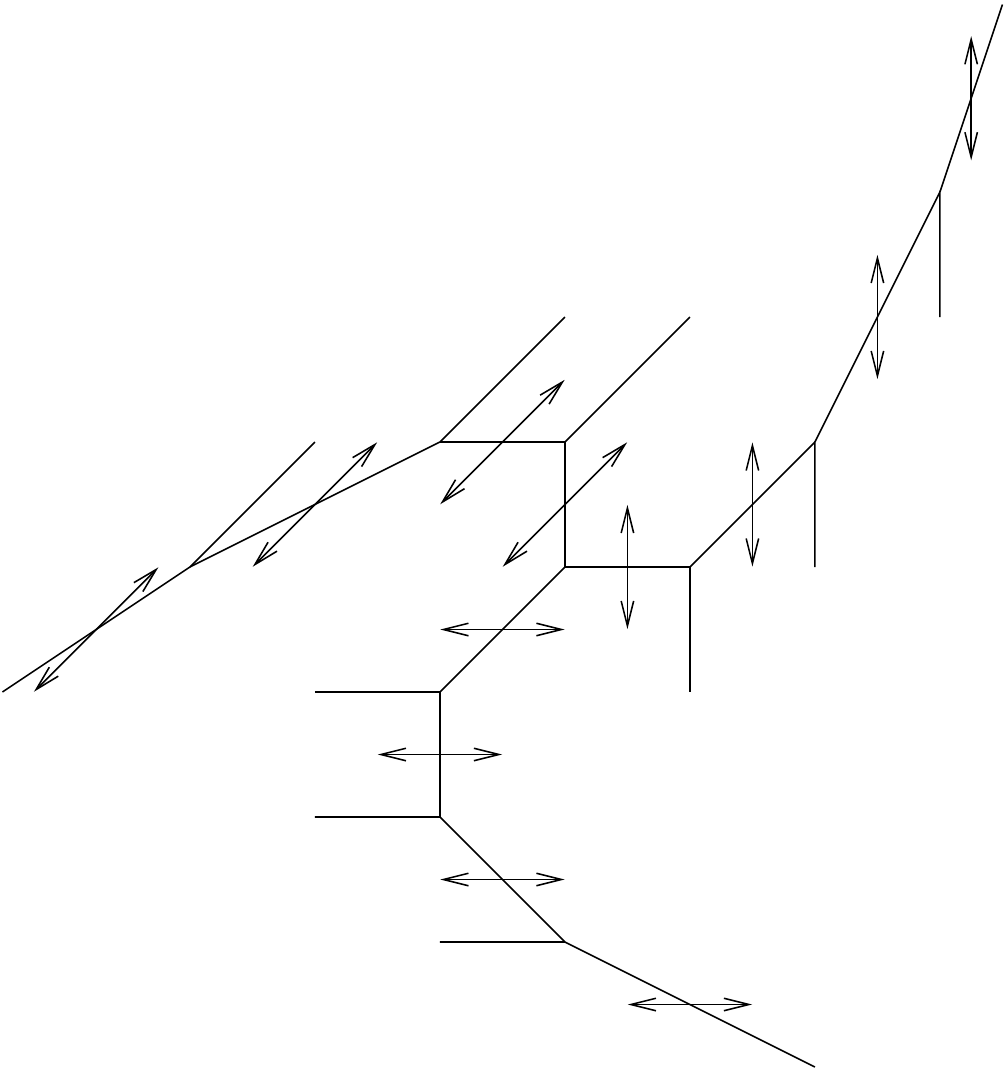}
\end{center}
\caption{Degeneration of $\hY_N$}
\end{figure}

We now degenerate each $\bP^1$ into two $\bP^1$'s intersecting
at a node. The total space of the normal bundle $\cO(n)\oplus \cO(-n-2)$
degenerates to $\cO(a)\oplus \cO(-a-1)$ and $\cO(b)\oplus \cO(-b-1)$ with $a+b=n$.
For each node we introduce a pair of framing vectors to encode the $\bS$--action
(see \fullref{fig:break}). We refer to \cite[Section 4]{LLLZ} for details.

The framing here corresponds to the framing of Lagrangian submanifolds
in the article by Aganagic, Klemm, Mari\~no and Vafa \cite{AKMV} and the framing of knots
and links in Chern--Simons theory.  In \fullref{fig:minimal}, all the $\bP^1$'s have
normal bundles $\cO(-1)\oplus \cO(-1)$ or $\cO\oplus \cO(-2)$; in Figure
15, all the $\bP^1$'s have normal bundles $\cO\oplus \cO(-1)$.

Using the connected version of the gluing formula \cite[Theorem 7.5]{LLLZ},
we have
\begin{multline*}
F^4_N(\la;\vt)
=\sum_{|\mu^i|>0}F_\vmu(\la;0,0,0)
\prod_{i=1}^3 \prod_{j=1}^{\ell(\mu^i)}
\sum_{k=1}^{N} e^{-\mu^i_j(t_{i,1}+\cdots+t_{i,k}) }\\
\cdot \left(\mu^i_jF_{(\mu^i_j), (\mu^i_j),\emptyset}(-1,0,0)\right)^{k-1}
\mu^i_jF_{(\mu^i_j),\emptyset,\emptyset}(\la;0,0,0)
\end{multline*}
$$= \sum_{|\mu^i|>0}F_\vmu(\la;0,0,0)
\prod_{i=1}^3\prod_{j=1}^{\ell(\mu^i)}\Bigl(\mu^i_j
\frac{(-1)^{\mu^i_j-1}\sqrt{-1}}{\mu^i_j[\mu^i_j]}
\sum_{k=1}^N e^{-\mu^i_j(t_{i,1}+\cdots +
  t_{i,k})}\Bigr)$$
So
\begin{multline}\label{eqn:fourI}
F^4_N(\la;\vt)=\\
\sum_{|\mu^i|>0} F_{\vmu}(\la;0,0,0)
\prod_{i=1}^3 (-1)^{|\mu^i|}(-\sqrt{-1})^{\ell(\mu^i)}
\prod_{j=1}^{\ell(\mu^i)}\biggl(
\frac{1}{[\mu^i_j]} \sum_{k=1}^N e^{-\mu^i_j(t_{i,1}+\cdots +t_{i,k}) } 
\biggr)
\end{multline}
In particular, using $F_{(1),(1),(1)}(\lam;0,0,0)=-\sqrt{-1}[1]$ and
\eqref{eqn:fourI}, we can recover \fullref{thm:tri}. Equation
\eqref{eqn:fourI} is equivalent to 
\begin{equation}\label{eqn:four}
F^4_N(\la;\vt)=\sum_{|\mu^i|>0}\tF_{\vmu}(\la;0,0,0) 
\prod_{i=1}^3 (-1)^{\ell(\mu^i)}e^{-|\mu^i| t_{i,1}}u^i_{\mu^i}(q,\vt_i)
\end{equation}
where $u^i_\mu(q;\vt_i)$'s are defined as in \fullref{sec:notation}, and
$$
\tF_{\vmu}(\la;0,0,0)=(-1)^{\sum_{i=1}^3|\mu^i|}\sqrt{-1}^{\ell(\vmu)} F_{\vmu}(\la;0,0,0).
$$
Similarly,
\begin{eqnarray}
\label{eqn:three}
\qquad F^3_N(\la;\vt_1,\vt_2)&=&\!\!\!\!
  \sum_{|\mu^1|>0,|\mu^2|>0}\!\!\!\!
  \tF_{\mu^1,\mu^2,\emptyset}(\la;0,0,0)
  \prod_{i=1}^2 (-1)^{\ell(\mu^i)}e^{-|\mu^i| t_{i,1} }
  u^i_{\mu^i}(q;\vt_i).\\
\notag
F^1_N(\la;\vt_1)&=&\sum_{\mu\neq \emptyset}
  \tF_{\mu,\emptyset,\emptyset} (\la;0,0,0)
  (-1)^{\ell(\mu)}e^{-|\mu| t_{1,1} }u^1_\mu(q;\vt_1)\\
\notag
&=&\sum_{\mu\neq \emptyset}F_{\mu,\emptyset,\emptyset}(\la;0,0,0)
  (-1)^{|\mu|}(-\sqrt{-1})^{\ell(\mu)}e^{-|\mu| t_{1,1} }u^1_{\mu}(q;\vt_1)\\
\notag
&=&\sum_{n>0} \frac{(-1)^{n-1}\sqrt{-1}}{n[n]} (-1)^n(-\sqrt{-1})
  e^{-nt_{1,1}} u^1_n(q;\vt_1).
\end{eqnarray}
So
\begin{equation}\label{eqn:one}
F^1_N(\la;\vt_1)= \sum_{n>0}
\frac{-1}{n[n]^2}
\sum_{k=1}^N e^{-n(t_{1,1}+\cdots +t_{1,k})}.
\end{equation}
This proves \fullref{thm:Fone}.

\begin{eqnarray*}
F^2_N(\la;\vt_1)&=&\sum_{n>0}F_{(n),\emptyset,\emptyset}(\la;-1,0,0)
\frac{(-1)^{n-1}\sqrt{-1}}{[n]}
\sum_{2\leq k_1\leq k_2\leq N} e^{-n(t_{1,k_1}+\cdots +t_{1,k_2})} \\
&=& \sum_{n>0}\frac{(-1)^n\sqrt{-1}}{n[n]}\cdot
\frac{(-1)^{n-1}\sqrt{-1}}{[n]}
\sum_{2\leq k_1\leq k_2\leq N} e^{-n(t_{1,k_1}+\cdots +t_{1,k_2})} 
\end{eqnarray*}
So
\begin{equation}\label{eqn:two}
F^2_N(\la;\vt_1)= \sum_{n=1}^\infty
\frac{1}{n[n]^2}
\sum_{2\leq k_1\leq k_2\leq N} e^{-n(t_{1,k_1}+\cdots +t_{1,k_2})} 
\end{equation}
This proves \fullref{thm:Ftwo}.

From \eqref{eqn:four}, \eqref{eqn:three}, and \eqref{eqn:one}, it
is clear that if $d_{i,1}>0$, then $\tN_{\vd}^g(\hY^N)=0$ unless
$$
d_{i,1}\geq d_{i,2}\geq \ldots\geq d_{i,k_i}>0.
$$
This proves \fullref{thm:decreasing}.

\begin{figure}[h]\label{fig:snake}
\begin{center}
\includegraphics[scale=0.77]{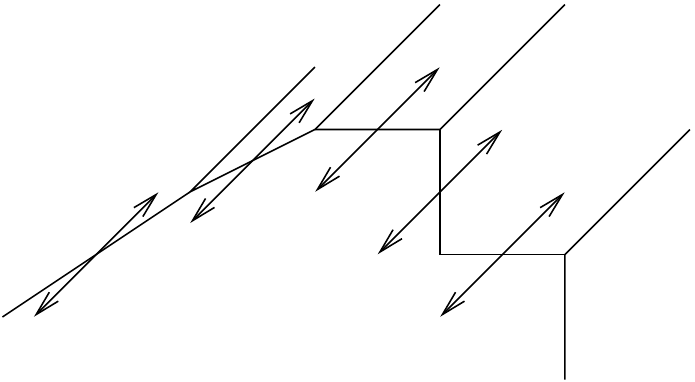}
\end{center}
\caption{Another degeneration}
\end{figure}

Let $F^5_N(\la;t_{1,1},\vt_2)$ be the contribution to 
$F^3_N(\la;\vt_1,\vt_2)$
from the case in \fullref{thm:Fthree}. To compute
$F^5_N(\la;t_{1,1},\vt_2)$, we consider the degeneration in 
\fullref{fig:snake}.

We have
\begin{multline*}
F^5_N(\la;t_{1,1},\vt_2)=\sum_{n>0}
F_{(n),\emptyset,\emptyset}(-1,0,0) n
F_{(n),(n),\emptyset}(-1,0,0)\\
\cdot\biggl(\sum_{k=1}^N (n F_{(n),(n),\emptyset}(-1,0,0 ) )^{k-1}
n F_{(n),\emptyset,\emptyset}(0,0,0)
e^{-nt_{1,1}-n(t_{2,1}+\cdots +t_{2,k}) }\biggr)\\
=\sum_{n>0}\frac{(-1)^n\sqrt{-1}}{n[n]}\cdot \frac{(-1)^{n-1}\sqrt{-1}}{[n]}
\biggl(\sum_{k=1}^N e^{-nt_{1,1}-n(t_{2,1}+\cdots +t_{2,k})}\biggr)
\end{multline*}

So 
\begin{equation}\label{eqn:five}
F^5_N(\la;t_{1,1},\vt_2)=\sum_{n=1}^\infty \frac{1}{n[n]^2}
\sum_{k=1}^N e^{-nt_{1,1}-n(t_{2,1}+\cdots +t_{2,k}) }
\end{equation}
\fullref{thm:Fthree} follows from \eqref{eqn:five}.

We now prove \fullref{thm:Ffour}.
Suppose that $d_{i,j}=d_i>0$ for $j\leq k_i$ and $d_{i,j}=0$ for $j>k_i$,
where $i=1,2,3$ and $1\leq k_i\leq N$. From \eqref{eqn:three} and
\eqref{eqn:four} it is easy to see that
\begin{eqnarray*}
\tN^g_{\vd_1,\vd_2,0}(\hY^N)&=&\tN^g_{d_1,d_2,0}(\hY^1),\\
\tN^g_{\vd_1,\vd_2,\vd_3}(\hY^N)&=&\tN^g_{d_1,d_2,d_3}(\hY^1).
\end{eqnarray*}
Recall that
$$
\tN^g_{\lo{d}}(\hY^1)=N^g_{\lo{d}}(Y^1)
$$
where $N^g_{\lo{d}}(Y^1)$ are given by \fullref{thm:threePone}.
\fullref{thm:Ffour} follows from \fullref{thm:Fthree}
and \fullref{thm:threePone}.

Finally, we prove \fullref{thm:tZ}. Let 
$$
F_N'(\la;\vt)=F_N(\la;\vt)-\sum_{i=1}^3 F_N^2(\la;\vt_i).
$$
Then
$$
F_N'(\la;\vt)=\sum_{\vmu \neq \vec{\emptyset}}\tF_{\vmu}(\la;0,0,0)
\prod_{i=1}^3 (-1)^{\ell(\mu^i)}e^{-|\mu^i| t_{i,1} } u^i_{\mu^i}(q;\vt^i).
$$
It was proved in \cite{LLLZ} that
\begin{equation}\label{eqn:tFW}
\tF^\bu_{\vmu}(\la;0,0,0)=\sum_{|\nu^i|=|\mu^i|}\tilde{\cW}_\vnu(q) 
\prod_{i=1}^3\frac{\chi_{\mu^i}(\nu^i)}{z_{\mu^i}}
\end{equation}
So 
\begin{eqnarray*}
\exp(F_N'(\la;\vt))&=&\sum_{\vnu}\tilde{\cW}_{(\nu^1)^t,(\nu^2)^t,(\nu^3)^t}(q)
\prod_{i=1}^3 (-1)^{|\nu^i|} e^{-|\nu^i| t_{i,1} } s_{\nu^i}(u^i(q,\vt_i))\\
&=&\sum_{\vmu}\tilde{\cW}_\vmu(q)
\prod_{i=1}^3 (-1)^{|\mu^i|} e^{-|\mu^i| t_{i,1} } s_{(\mu^i)^t}(u^i(q,\vt_i))
\end{eqnarray*}
where $s_\mu(u^i(q,\vt_i))$'s are defined as in \fullref{sec:notation}.
We conclude that
\begin{multline*}
\exp(F_N(\la;\vt))=\\
\exp\Bigl(\sum_{i=1}^3 F_N^2(\lam;\vt_i)\Bigr)
\sum_\vmu \tilde{\cW}_{\vmu}(q) 
\prod_{i=1}^3(-1)^{|\mu^i|} e^{-|\mu^i|t_{i,1}} s_{(\mu^i)^t}(u^i(q,\vt_i)).
\end{multline*}
This completes the proof of \fullref{thm:tZ}.

\section{Physical theory of the topological vertex}\label{sec:vertex}
In this section, we compute the local Gromov--Witten invariants considered in this paper 
by using the physical theory of the topological vertex. Typical computations in this theory 
involve formal sums over Young tableaux, and in some cases, like the one considered here, 
it is more convenient to use the operator formalism on Fock spaces. After a short overview of this formalism, we will 
use it to compute the partition functions for local Gromov--Witten invariants.  

\subsection{Operator formalism}
We introduce:
$$
|p_\mu\rangle = \prod_{i=1}^{\ell(\mu)} \alpha_{-\mu_i} |0\rangle,\ \
|s_\mu\rangle = \sum_{|\nu|=|\mu|} \frac{\chi_\mu(\nu)}{z_\nu}
|p_\nu\rangle,\ \
|s_{\lambda/\mu}\rangle=\sum_\nu c^\lambda_{\mu\nu} |s_\nu\rangle
$$
where $\alpha_n$ satisfy the commutation relations
$$[\alpha_m,\alpha_n]=m\delta_{m+n,0},$$
and for $n>0$, $\alpha_n |0\rangle=0$. The dual vector space is obtained by acting with the
operators $\alpha_n$ on the state $\langle 0|$, and the pairing
is defined by $\langle 0 | 0 \rangle=1$. One then finds,
$$\langle p_\mu| p_\nu\rangle= z_\mu \delta_{\mu,\nu},\quad
  \langle s_\mu|s_\nu \rangle= \delta_{\mu,\nu}.$$
The coherent state $|t \rangle$ is defined as
$$|t\rangle =\exp\biggl(\sum_{n=1}^\infty \frac{t_n}{n}\alpha_{-n}\biggr)
  |0\rangle =\sum_\mu\frac{t_\mu}{z_\mu} |p_\mu\rangle$$
where
$$t_\mu= t_{\mu_1}\cdots t_{\mu_{\ell(\mu)}}$$
and one has
$$\langle s | t \rangle = \exp\biggl(\sum_{n=1}^\infty
  \frac{s_n t_n}{n}\biggr).$$
The elements $|s_{\mu}\rangle$, where $\mu$ is a partition,
span a vector space ${\cal H}$ that can be
identified with the ring of symmetric functions $\Lambda$ in an infinite number of variables,
and therefore it inherits a ring structure from $\Lambda$.
This identification can be made by considering the map
$$|s_{\mu}\rangle \rightarrow \langle t | s_{\mu}\rangle \equiv s_{\mu}(t),$$
where $s_{\mu}(t)$ gives the Schur function after identifying
$t_n=p_n=x_1^n + x_2^n + \cdots$.

Given a coherent state $|t \rangle$, it is useful to define the coherent states
$$|t^{\omega} \rangle, \qquad |t^{\zeta} \rangle, \qquad |t^{\xi} \rangle$$
where
$$t^{\omega}_n=(-1)^{n+1} t_n, \qquad t^{\zeta}_n=-t_n, \qquad
  t^{\xi}_n=(-1)^n t_n. $$
Using $\chi_{\mu^t}(\nu)=(-1)^{|\nu| + \ell(\nu)} \chi_{\mu}(\nu)$ it is easy to
show that
\begin{equation}
s_{\mu}(t^{\omega})=s_{\mu^t}(t), \qquad s_{\mu}(t^{\xi})=(-1)^{|\mu|}s_{\mu}(t).
\label{transposition}
\end{equation}
The ring of symmetric polynomials is endowed with a coproduct
structure (see for example Macdonald \cite[Ex. 25 of I.5]{M})
$$\Delta\co  \Lambda  \rightarrow \Lambda\otimes \Lambda$$
which is a ring homomorphism, and is defined by
$$\Delta(s_{\lambda})=\sum_{\mu} s_{\lambda/\mu}\otimes s_{\mu}.$$
The $n$th power sums $p_n$ are primitive
elements  of $\Lambda$ under this coproduct, and one has
$$\Delta p_n=p_n \otimes 1 + 1 \otimes p_n.$$
We then have a inherited coproduct $\Delta\co  {\cal H} \rightarrow
{\cal H}\otimes {\cal H}$, and
it is easy to see that it acts as follows on coherent states:
$$\Delta(|t \rangle)=|t \rangle \otimes |t \rangle.$$
This gives the following identity, which will be 
useful in proving \fullref{thm:Zone}:
\begin{equation}
\label{coproduct}
\sum_{\mu, \nu} s_{\mu}(t) | s_{\mu/\nu}\rangle\langle s_{\nu}|=|t \rangle \langle t |.
\end{equation}
We need now explicit expressions for ${\cal W}_{\mu} (q)$ and ${\cal W}_{\mu \nu}(q)$ in the
operator formalism. Using \eqref{qdimension}, one immediately finds
\begin{equation}
{\cal W}_{\mu} (q)=s_{\mu}(\beta), \quad 
\beta_n=\frac{1}{[n]}=\frac{1}{q^{n/2} - q^{-n/2} },
\label{unknot}
\end{equation}
therefore
\begin{equation}
\sum_{\mu}{\cal W}_{\mu} (q) |s_{\mu}\rangle=|\beta\rangle,
\label{diskstate}
\end{equation}
where $|\beta \rangle$ is a coherent state with  $
\beta_n$ given in \eqref{unknot}.
One can then write
$${\cal W}_{\mu}(q)=\langle s_{\mu}| \beta\rangle.$$
We introduce the operator $q^{\pm \kappa/2}$ defined by
$$q^{\pm \kappa/2}|s_{\mu}\rangle= q^{\pm \kappa_{\mu}/2} |s_{\mu}\rangle.$$
We also define an operator $W$ as
$${\cal W}_{\mu\nu}(q)=\langle s_\mu | W | s_\nu \rangle.$$
In the proof of \fullref{thm:Zone}, we will need an explicit expression for
$$\langle t | q^{-\kappa/2} W q^{-\kappa/2} |{\bar t}\rangle,$$
where $|t\rangle$, $|\bar t\rangle$ are coherent states. Using \eqref{symw} one finds
\begin{eqnarray*}
\langle t | q^{-\kappa/2} W q^{-\kappa/2} |{\bar t}\rangle
&=& \sum_{\mu, \nu} q^{-\kappa_{\mu}/2 - \kappa_{\nu}/2}{\cal W}_{\mu,\nu}(q) s_{\mu}(t)
s_{\nu}(\bar t)\\
&=&\sum_{\mu,\nu,\sigma} s_{\mu^t}(t^{\omega}) {\cal W}_{\mu^t/\sigma}(q)
{\cal W}_{\nu^t/\sigma}(q)s_{\nu^t}(\bar t^{\omega})\\
&=&\sum_{\mu,\nu,\sigma,\tau} s_{\mu}(t^{\omega}) \langle \beta|s_{\mu/\sigma}\rangle \langle
s_{\sigma} | s_{\tau}\rangle\langle s_{\nu/\tau}|\beta \rangle s_{\nu}(\bar t^{\omega})\\
&=&\langle \beta | t^{\omega}\rangle\langle  t^{\omega} | \bar t^{\omega}\rangle\langle \bar t^{\omega}
|\beta \rangle,
\end{eqnarray*}
where in the last step we have used \eqref{coproduct} twice. The last quantity is expressed solely
in terms of products of coherent states, and we finally find
\begin{equation}
\label{generating}
\langle t | q^{-\kappa/2} W q^{-\kappa/2} |{\bar t}\rangle=
\exp\biggl(\sum_{n=1}^\infty \frac{(-1)^{n+1}(t_n +{\bar t}_n) }{n [n]} + \frac{t_n \bar t_n}{n}
\biggr).
\end{equation}
This result was previously obtained, in slightly different form, by
Aganagic, Dijkgraaf, Klemm, Mari\~no and Vafa \cite{ADKMV}, and Zhou
\cite{ZI}. 
It is also possible to compute $\langle t| q^{-\kappa/2} W q^{-\kappa/2} |s_{\mu}\rangle$ by following the same steps. One finds 
$$\langle t| q^{-\kappa/2} W q^{-\kappa/2} |s_{\mu}\rangle=\langle \beta
  |t^{\omega}\rangle \sum_{\nu} {\cal W}_{\mu^t/\nu}(q) s_{\nu}(t^{\omega}).$$
The sum over the representation $\nu$ can be performed explicitly by using
the following formula (see Macdonald \cite{M}):
$$
\sum_{\nu} s_{\mu/\nu}(x) s_{\nu}(y)=s_{\mu}(x,y), 
$$
where $x$, $y$ are variables of the Schur polynomials. We then find
\be
\langle t| q^{-\kappa/2} W q^{-\kappa/2} |s_{\mu}\rangle
=\langle \beta | t^{\omega}\rangle s_{\mu}(u_n), \qquad 
u_n=\frac{(-1)^{n+1}}{[n]} + t_n.
\label{mixedsande}
\ee

\subsection{The closed topological vertex}\label{sec:NW}

The physical theory of the topological vertex \cite{AKMV}
gives the following expression for $\tZ_1(\la;\vt)$:
\begin{align*}
&Z_1(\la;\vt) = \\
&\!\sum_{\mu^1,\mu^2,\mu^3}
\!\!\!\!{\cal W}_{(\mu^3)^t\!,(\mu^2)^t\!,(\mu^1)^t}(q)
 {\cal W}_{\mu^1}(q) {\cal W}_{\mu^2}(q) {\cal W}_{\mu^3}(q)
     e^{-\!\sum_{i=1}^3|\mu^i| t_i}
     (-1)^{\sum_{i=1}^3|\mu^i|} = \\
&\!\sum_{\mu^1,\mu^2,\mu^3}
\!\!\!\!{\cal W}_{\vec \mu}(q)
   {\cal W}_{(\mu^1)^t}(q) {\cal W}_{(\mu^2)^t}(q) {\cal W}_{(\mu^3)^t}(q)
     e^{-|\mu^1| t_1-|\mu^2| t_2-|\mu^3| t_3}
     (-1)^{|\mu^1|+|\mu^2|+|\mu^3|}
\end{align*}
where $\vt=(\lo{t})$, and in the last step we have used
$$
\cW_{(\mu^3)^t, (\mu^2)^t, (\mu^1)^t}(q)
= q^{-\sum_{i=1}^3 \kappa_{\mu^i}/2}\cW_{\up{\mu}}(q), \quad
\cW_{\mu}(q)= q^{\kappa_\mu/2}\cW_{\mu^t}(q).
$$
In this subsection, we will prove \eqref{eqn:Zone} in
\fullref{sec:introduction}:
\begin{prop}\label{thm:Zone}
$$
Z_1(\la;\vt)=\exp\biggl(\sum_{n=1}^\infty \frac{Q_n(\vt)}{-n[n]^2}\biggr)
$$
where
$$Q_n(\vt)=e^{-nt_1}+ e^{-nt_2}+ e^{-nt_3}
  -e^{-n(t_1+t_2)}-e^{-n(t_2+t_3)} -e^{-n(t_3+t_1)}
  +e^{-n(t_1+t_2+t_3)}.$$
\end{prop}

\begin{proof}
\begin{align}
Z_1&(\la;\vt) = \nonumber \\
&\sum_{\scriptsize\begin{array}{c}
\mu^1\!,\mu^2\!,\mu^3\!,\\\rho,\rho^1\!,\rho^3\!,\rho'
\end{array}} \!\!\!\!  \!\!\!\!
c^{\mu^1}_{\rho \rho^1}  {\cal W}_{(\mu^1)^t}(q)
\bigl\langle s_{(\mu^2)^t} \big| W\big| s_{\rho^1}\bigr\rangle
\bigl\langle s_\rho\big| {\cal W}_{(\mu^3)^t}(q)
  c^{\mu^3}_{\rho'\rho^3}\big| s_{\rho'}\bigr\rangle
\bigl\langle s_{\rho^3} \big| W \big| s_{\mu^2}\bigr\rangle\nonumber \\[-15pt]
&\hskip 100pt\cdot e^{-|\mu^1| t_1} (-1)^{|\mu^1|}
\cdot e^{-|\mu^2| t_2} (-1)^{|\mu^2|}
\cdot e^{-|\mu^3| t_3} (-1)^{|\mu^3|} = \nonumber \\
&\sum_{\mu^2}(-1)^{|\mu^2|} e^{-|\mu^2|t_2}
\Bigl\langle s_{(\mu^2)^t}\Big| W
\Bigl( \!\!\sum_{\rho,\rho^1\!,\mu^1}\!\!\!\! c^{\mu^1}_{\rho,\rho^1} (-1)^{|\mu^1|} e^{-|\mu^1|t_1}
{\cal W}_{(\mu^1)^t}(q)| s_{\rho^1}\rangle \langle s_\rho |\Bigr)
\label{zqfirst} \\
&\hskip 80pt\Bigl(\!\!\sum_{\rho'\!,\rho^3\!,\mu^3} \!\! c^{\mu^3}_{\rho'\rho^3}
(-1)^{|\mu^3|} e^{-|\mu^3|t_3}{\cal W}_{(\mu^3)^t}(q)| s_{\rho'}\rangle
\langle s_{\rho^3}|
\Bigr) W\Big| s_{\mu^2}\Bigr\rangle.\nonumber
\end{align}
Here we have used the explicit expression for the topological vertex  \eqref{topvertex} and 
the identity \eqref{qdimtrans}. We now write
\begin{align}
\sum_{\rho,\rho^1,\mu}& c^{\mu}_{\rho,\rho^1} (-1)^{|\mu|} e^{-|\mu|t_1}
{\cal W}_{\mu^t}(q)| s_{\rho^1}\rangle \langle s_\rho | \nonumber \\
&=  \sum_{\rho,\mu} (-1)^{|\mu|} e^{-|\mu|t_1}
{\cal W}_{\mu^t}(q)| s_{\mu/\rho}\rangle \langle s_\rho | \nonumber \\
&= \sum_{\rho,\mu} s_{\mu}(u^\zeta)| s_{\mu/\rho}\rangle \langle s_\rho |
\nonumber \\
&= |u^{\zeta}\rangle\langle u^{\zeta}|,
\label{manip}
\end{align}
where
\begin{equation}
u_n=e^{-nt_1}\beta_n.
\end{equation}
and in the last step of \eqref{manip} use has been made of \eqref{coproduct} and of the fact that, under
$t_n \rightarrow a t_n$ one has
$$s_{\mu}(a t)=a^{|\mu|}s_{\mu}(t).$$
Following the same steps for the second bracket in \eqref{zqfirst} one finds,
$$Z_1(\la;\vt)=\sum_{\mu}(-1)^{|\mu|} e^{-|\mu|t_2}
  \langle s_{\mu^t}| W| u^\zeta\rangle
  \langle u^{\zeta}|v^{\zeta}\rangle\langle v^{\zeta}|W | s_{\mu}\rangle.$$
with
\begin{equation}
v_n=e^{-n t_3} \beta_n.
\end{equation}
Notice that
$$
|u^{\zeta}\rangle=\sum_{\mu} (-1)^{|\mu|}e^{-|\mu| t_1} {\cal W}_{\mu^t}(q) |s_{\mu}\rangle,
$$
and using again that ${\cal W}_{\mu^t}(q)=q^{-\kappa_{\mu}/2}{\cal W}_{\mu}(q)$ one can write
$$
|u^{\zeta}\rangle=q^{-\kappa/2} |u^{\xi}\rangle.
$$
We have similar equations for $|v^{\zeta}\rangle$. Since $\langle u^{\zeta}|v^{\zeta}\rangle$
is a product of coherent states, we only have to evaluate
$$
\sum_{\mu}(-1)^{|\mu|} e^{-|\mu|t_2}\langle v^{\xi}| q^{-\kappa/2}W q^{-\kappa/2}| s_{\mu}\rangle
\langle s_{\mu^t}| q^{-\kappa/2}Wq^{-\kappa/2} | u^{\xi} \rangle,
$$
where we used that $\kappa_{\mu}=-\kappa_{\mu^t}$. The last step involves writing
\begin{equation}
\sum_{\mu}(-1)^{|\mu|} e^{-|\mu|t_2} |s_{\mu} \rangle \otimes |s_{\mu^t} \rangle
=\exp\left(-\sum_{n=1}^\infty \frac{e^{-nt_2}}{n} \alpha^{(1)}_{-n}\otimes \alpha^{(2)}_{-n}
\right) |0\rangle_1 \otimes |0 \rangle_2
\label{finalstate}
\end{equation}
which is an element of ${\cal H}_1\otimes {\cal H}_2$, and we have introduced explicit indices
$1,2$ to label the factors in the tensor product. We first take the scalar product of this state
with $\langle v^{\xi}| q^{-\kappa/2}W q^{-\kappa/2} \in {\cal H}_1^{*}$ to obtain a state in
${\cal H}_2$. In order to do that, we can regard \eqref{finalstate} as a coherent state with
$t_n=-e^{-nt_2}\smash{\alpha^{(2)}_{-n}}$, therefore we can use the formula \eqref{generating} to
obtain the element in ${\cal H}_2$
\begin{eqnarray*}
&&
\exp\left(\sum_{n=1}^\infty -{v_n\over n[n]} +{(-1)^n e^{-nt_2}\over n[n]}\alpha^{(2)}_{-n}
+ {(-1)^{n+1}e^{-nt_2} v_n \over n} \alpha^{(2)}_{-n}
\right)|0 \rangle\\
&=&\exp\left(-\sum_{n=1}^\infty {e^{-nt_3}\over n[n]^2}\right) |w \rangle
\end{eqnarray*}
where $|w\rangle$ is a coherent state in ${\cal H}_2$ given by
$$
w_n=(-1)^{n+1} \beta_n e^{-nt_2}(e^{-nt_3}-1).
$$
The remaining step is to compute $\langle w| q^{-\kappa/2}W q^{-\kappa/2} | u^{\xi}\rangle$,
which can be done again with the help of \eqref{generating}. Collecting all terms, one finds
\begin{equation}
\begin{aligned}
& Z_1(\la;\vt)= \exp\biggl\{ -\sum_{n=1}^{\infty} {1\over n[n]^2}
(e^{-nt_1}+ e^{-nt_2} + e^{-n t_3}- e^{-n(t_1+t_2)}  \\
& \makebox[3.4cm]{ }
 -e^{-n(t_2+t_3)} -e^{-n(t_3+t_1)} + e^{-n(t_1+t_2+t_3)} )\biggr\}.
\end{aligned}
\end{equation}
This completes the proof.
\end{proof}

\subsection{A chain of rational curves}\label{sec:minustwo}

Let $F_N^2(\lambda;\vt_1)$ be defined as in \fullref{sec:cases},
where $N\geq 2$. It can  be viewed as a generating function of formal Gromov--Witten
invariants of a chain of $(N-1)$ rational curves with normal 
bundles $\cO\oplus \cO(-2)$.
In this section, we will compute
$$
Z_N^2(q;t_{1,2},\ldots,t_{1,N})=\exp(F_N^2(\lambda;\vt_1))
$$
using vertex techniques.

We have
\begin{multline}
Z_N^2(q;t_{1,2},\ldots,t_{1,N}) = \\
\sum_{\mu^2,\ldots,\mu^N} \prod_{i=2}^{N+1} {\cal W}_{(\mu^{i-1})^t, \emptyset , \mu^i}(q) 
    q^{-\kappa_{\mu^i}/2} e^{-|\mu^i|t_{1,i}} = \\
\sum_{\mu^2,\ldots,\mu^N} \prod_{i=2}^{N+1} q^{-\kappa_{\mu^{i-1}}/2} {\cal W}_{\mu^{i-1},\mu^i}(q) 
q^{-\kappa_{\mu^i}/2} e^{-|\mu^i|t_{1,i}},
\label{chainpf}
\end{multline}
where $\mu^1=\mu^{N+1}=\emptyset$. 
Using the above techniques, in particular \eqref{diskstate} and \eqref{transposition}, we
can write the above expression as
$$
Z^2_N(q;t_{1,2},\ldots,t_{1,N})
=\langle u^2| \biggl( \prod_{i=3}^N q^{-\kappa/2} W q^{-\kappa/2} {\cal O}_i \biggr) |u^{N+1}\rangle,
$$
where $| u^2\rangle$, $|u^{N+1}\rangle$ are coherent states defined by
$$
u^2_n={(-1)^{n+1} e^{-nt_{1,2} }\over [n]}, \quad 
u^{N+1}_n={(-1)^{n+1} \over [n] }
$$
and
$$
{\cal O}_i=\sum_{\mu^i} | s_{\mu^i}\rangle e^{-|\mu^i|t_{1,i} } \langle s_{\mu^i}|, \quad i=3, \cdots, N.
$$
We first compute
$$
\langle u^2 | q^{-\kappa/2} W q^{-\kappa/2} {\cal O}_2,
$$
which is an element of ${\cal H}^*$. We proceed as in the computation
following \eqref{finalstate}
above, to obtain
$$
\langle u^2 | q^{-\kappa/2} W q^{-\kappa_2} {\cal O}_3 = 
\exp\biggl( \sum_{n=1}^{\infty} {(-1)^{n+1} \over n[n]} u^2_n\biggr)
\langle u^3|,
$$
where $\langle u^3|$ is a coherent state defined by
$$
u^3_n=e^{-n t_{1,3}} \biggl( u_n^2 + {(-1)^{n+1} \over [n]}\biggr).
$$
We can now compute $Z_N^2(q;t_{1,2},\ldots,t_{1,N})$ recursively, 
defining the coherent state $|u^i \rangle$ as
$$
\langle u^{i-1} | q^{-\kappa/2} W q^{-\kappa/2} {\cal O}_i = 
\exp\biggl( \sum_{n=1}^{\infty} {(-1)^{n+1} \over n[n]}u^{i-1}_n\biggr)
\langle u^i|,
$$
where
 \begin{equation}
u^i_n=e^{-n t_{1,i}} \biggl( u_n^{i-1} + {(-1)^{n+1} \over [n]}\biggr), \quad i=3, \cdots, N.
\label{recursivedef}
\end{equation}
One then finds, by using \eqref{generating} repeatedly, that
\begin{eqnarray*}
Z^2_N(q;t_{1,2},\ldots,t_{1,N})
&=&\exp \biggl( \sum_{n=1}^{\infty} {(-1)^{n+1}\over n[n] }\sum_{i=2}^{N-1} u^i_n +
\sum_{n=1}^\infty {1\over n}u^N_n u^{N+1}_n
\biggr)\\
&=&\exp \biggl( \sum_{n=1}^{\infty} {(-1)^{n+1}\over n[n] }\sum_{i=2}^{N} u^i_n \biggr).
\end{eqnarray*}
The recursion relation defining $u^i_n$ is easily solved:
$$
u^i_n=\frac{(-1)^{n+1}}{[n]}
\Bigl(\sum_{j=2}^i e^{-n(t_{1,j}+\cdots +t_{1,i})}\Bigr) , \quad i=2, \cdots, N,
$$
and putting everything together we finally obtain
$$
Z^2_N(q;t_{1,2},\ldots,t_{1,N})=
\exp\biggl\{ \sum_{n=1}^{\infty} {1\over n [n]^2} 
\Bigl( \sum_{2\le i \le j\le N}
e^{-n(t_{1,i}+\cdots+t_{1,j})}\Bigr) \biggr\},
$$
or equivalently,
$$
F^2_N(\lambda;\vt_1)=
\sum_{n=1}^{\infty} {1\over n [n]^2} 
\Bigl( \sum_{2\le i \le j\le N}
e^{-n(t_{1,i}+\cdots+t_{1,j})}\Bigr)
$$
which agrees with \eqref{eqn:Ftwo}.

Notice that the non-trivial Gopakumar--Vafa invariants for this geometry occur for K\"ahler
classes which are in one-to-one correspondence with the positive roots of the Lie algebra $A_{N-1}$.

\subsection{Minimal trivalent configuration}

Let us finally consider the minimal trivalent configuration. 
We will allow the three chains of $\bP^1$'s to have different
lengths $\lo{N}$:
\[
Y^{\lo{N}}=\bigcup_{1\leq i\leq N_1} A_i\cup
           \bigcup_{1\leq j\leq N_2} B_j \cup
           \bigcup_{1\leq k\leq N_3} C_j.
\]
So we have
\[
\begin{array}{cc}
\vd=(\lo{\vd}), & \vd_i=(d_{i,1}\ldots,d_{i,N_i}),\\
\vt=(\lo{\vt}), & \vt_i=(t_{i,1},\ldots,t_{i,N_i}),
\end{array}
\]
and we define
$$
u^i_n(q;\vt)=\frac{1}{[n]}\Bigl( 1+\sum_{k=2}^{N_i} e^{-n(t_{i,2}+\cdots+t_{i,k})} \Bigr).
$$

The rules of the topological vertex give the following expression, 
\begin{multline*}
 Z_{N_1,N_2,N_3}(\la;\vt)\\
=\sum_{\mu^{i,j}} (-1)^{\sum_{i=1}^3 
|\mu^{i,1}|} e^{-\sum_{i=1}^3 |\mu^{i,1}||t_{i,1}} q^{\sum_{i=1}^3 \kappa_{\mu^{i,1}}/2} 
{\cal W}_{(\mu^{3,1})^t, (\mu^{2,1})^t, (\mu^{1,1})^t}(q) \\[-10pt]
  \cdot \prod_{i=1}^3 \prod_{j=2}^{N_i+1} 
q^{-\kappa_{\mu^{i,j-1}}/2} {\cal W}_{\mu^{i,j-1},\mu^{i,j}}(q) q^{-\kappa_{\mu^{i,j}}/2}(q) e^{-|\mu^{i,j}|t_{i,j}}
\end{multline*}
where $\mu^{i, N_i+1}=\emptyset$. We will show that this expression can
be simplified as follows:
\begin{prop}\label{thm:Z}
\begin{multline}\label{eqn:Nthree}
Z_{\lo{N}}(\la;\vt)=\exp\biggl(\sum_{n=1}^\infty
  \frac{1}{n[n]^2}\sum_{i=1}^3\sum_{2\leq k_1\leq k_2 \leq N_i}
  e^{-n(t_{i,k_1}+\cdots +t_{i,k_2}) }\biggr) \\
\cdot \sum_{\vmu}\cW_{\mu^1,\mu^2,\mu^3}(q) 
  \prod_{i=1}^3 (-1)^{|\mu^i|} e^{-|\mu^i|t_{i,1}}s_{(\mu^i)^t}(u^i(q,\vt_i)))
\end{multline}
\end{prop}
Equation \eqref{eqn:Z} in \fullref{sec:introduction} corresponds to the case $N_1=N_2=N_3=N$.

\begin{proof}
The sum over the partitions $\mu^{i,j}$, $2\leq j\leq N_i$, 
can be performed by following the same steps that we made before, and
making use of \eqref{mixedsande}. 
After writing $\mu^{i,1}\rightarrow (\mu^{i})^t$, the resulting expression takes the following form:
\begin{multline}
Z_{\lo{N}}(\la;\vt)= \exp\biggl\{ \sum_{n=1}^{\infty}
  {(-1)^{n+1}\over n[n]} \sum_{i=1}^3 \sum_{j=2}^{N_i} v_n^{i,j} \biggr\}\\
\cdot\sum_{\vmu} (-1)^{\sum_{i=1}^3 |\mu^i|} e^{-\sum_{i=1}^3 |\mu^i| t_{i,1} }
  {\cal W}_{\mu^3,\mu^2, \mu^1}(q) 
 \prod_{i=1}^3 q^{-\kappa_{\mu^i}/2} s_{\mu^i}(u^i).
 \label{mtc}
\end{multline}
 In \eqref{mtc}, the variables $v_n^{i,j}$ are defined recursively by
\begin{align*}
v_n^{i,N_i}&={(-1)^{n+1} e^{-nt_{i,N_i}} \over [n]},\\
v_n^{i,j-1}&=e^{-n t_{i,j-1}}
  \biggl( v_n^{i,j} + {(-1)^{n+1} \over [n]}\biggr), \quad j=3, \ldots, N_i, \\
u^i_n&= {1\over [n]} + (-1)^{n+1} v_n^{i,2} 
\end{align*}
The recursion defining $v^{i,j}_n$ can be easily solved:
\begin{align*}
v^{i,j}_n&=\frac{(-1)^{n+1}}{[n]}\Bigl( \sum_{k=j}^{N_i} e^{-n(t_{i,j}+\cdots+t_{i,k})}\Bigr),
\quad j=2,\ldots, N_i,\\
u^i_n&=\frac{1}{[n]}\Bigl( 1+\sum_{k=2}^{N_i} e^{-n(t_{i,2}+\cdots+t_{i,k})}\Bigr).
\end{align*}
Therefore
\begin{multline}
Z_{\lo{N}}(\la,\vt)= \exp\biggl\{ \sum_{n=1}^{\infty} {1\over n[n]^2} 
  \sum_{i=1}^3 \sum_{2\leq j\leq k\leq N_i}
  e^{-n(t_{i,j}+\cdots+ t_{i,k})} \biggr\} \\
\cdot\sum_{\vmu} (-1)^{\sum_{i=1}^3 |\mu^i|}e^{-\sum_{i=1}^3|\mu^i|t_{i,1}}
  {\cal W}_{\mu^3,\mu^2,\mu^1}(q) 
  \prod_{i=1}^3 q^{-\kappa_{\mu^i}/2} s_{\mu^i}(u^i).
\label{eqn:mtc}
\end{multline}
Recall that
$$
\cW_{(\mu^1)^t,(\mu^2)^t,(\mu^3)^t}(q)= q^{-\sum_{i=1}^3 \kappa_{\mu^i}/2} \cW_{\mu^3,\mu^2,\mu^1}(q),
$$
so \eqref{eqn:mtc} is equivalent to \eqref{eqn:Nthree}.
\end{proof}

The closed topological vertex (\fullref{sec:NW}) and
chain of rational curves (\fullref{sec:minustwo}) can be obtained
taking limits of $Z_{\lo{N}}(\la;\vt)$:
\begin{enumerate}
\item Let $t_{i,j} \rightarrow \infty$ for $i=1,2,3$, $j\geq 2$. One has 
$$
v_n^{i,j}=0,  \quad i=1,2,3,\quad j\ge 2.
$$
and $u_n^i= \beta_n$ for $i=1,2,3$,  so  
$s_{\mu^i}(u^i)= \cW_{\mu^i}(q)= q^{\kappa_{\mu^i}/2}\cW_{(\mu^i)^t}(q)$
and $Z_{\lo{N}}(\la;\vt)$ becomes the closed topological vertex.
\item   Let $t_{i,j}\rightarrow \infty$ for $i=2,3$, and $t_{1,1}\rightarrow \infty$.
        We recover a chain of spheres with K\"ahler parameters $t_{1,2}, \ldots, t_{1,N_1}$. 
\end{enumerate}

 Unfortunately, the sum over partitions in \eqref{eqn:mtc} 
 can not be evaluated in close form as we did before. In fact, explicit computations show that 
 $Z_{\lo{N}}(\la;\vt)$ involves Gopakumar--Vafa invariants at higher genera, 
 and seem to indicate that there are infinitely many degrees $d_{i,j}$ for which 
 the Gopakumar--Vafa invariants are non-vanishing. If we write  
$Z_{\lo{N}} (\la; \vt)$ in the Gopakumar--Vafa 
 form
 $$
 Z_{\lo{N}} (\la;\vt)=\exp\biggl(\sum_{\ell=1}^{\infty} 
\sum_{g=0}^{\infty} \sum_{{\bf d}} n^{g}_{\bf d} {1\over \ell}[\ell]^{2g-2} e^{-\ell {\bf d} \cdot {\vt}}\biggr)
 $$
 where  ${\bf d} \cdot {\vt}=\sum_{i,j} d_{i,j} t_{i,j}$, then, for $N=2$ we find for example
 $$
 \begin{aligned}
 n^0_{(1,1),(1,1),(1,1)}=&-1,\\
 n^0_{(2,1),(1,1),(1,1)}=&1,\\
 n^0_{(1,0), (2,1),(2,1)}=&-2,\\
 n^0_{(1,1),(2,1),(2,1)}=&-2, \\ 
 n^0_{(2,1),(2,1),(2,1)}=&4,
 \end{aligned}
 $$
and they vanish for $g>0$. Due to the cyclic symmetry of the configuration, the same values are obtained for cyclic permutations of the three sets of degrees.  

If, say, $t_{3,j} \rightarrow \infty$ for $j\ge 1$, so we have two lines of spheres joined by a 
 two-vertex, then one can perform the sum over one of the two remaining partitions. This is because 
 ${\cal W}_{\emptyset, \mu_{2},\mu_1}(q)={\cal W}_{\mu^2, (\mu^1)^t}(q)
 q^{\kappa_{\mu^1}/2}$, and \eqref{mtc} reads
\begin{multline*}
Z_{N_1,N_2}(\la;\vt_1, \vt_2)=
  \exp\biggl\{ \sum_{n=1}^{\infty} {(-1)^{n+1}\over n[n]} 
  \sum_{i=1}^2 \sum_{j=2}^{N_i} v_n^{i,j} \biggr\}\\
\sum_{\mu^1, \mu^2}  (-1)^{\sum_{i=1}^2 |\mu^i|}
  e^{-\sum_{i=1}^2|\mu^i|t_{i,1} } 
  s_{\mu^2} (u^2) q^{-\kappa_{\mu^2}/2}
  {\cal W}_{\mu^2, (\mu^1)^t}(q) s_{\mu^1}(u^1).
\end{multline*}
Using again \eqref{mixedsande}, and relabelling $\mu^1\rightarrow \mu$, we finally obtain
\begin{multline*}
Z_{N_1,N_2}(\la;\vt_1, \vt_2)=
  \exp \biggl( \sum_{n=1}^{\infty} {(-1)^{n+1}\over n[n] }\Bigl(v^{2,1}_n 
  +\sum_{i=1}^2 \sum_{j=2}^{N_i} v^{i,j}_n\Bigr)
  +{1\over n}v^{2,1}_n v^{2,2}_n \biggr) \\
\sum_{\mu} q^{-\kappa_{\mu}/2} (-1)^{|\mu|} e^{-t_{1,1} |\mu|}
  s_{\mu}(\hat u^2)s_{\mu}(u^1),
\end{multline*}
where 
$$
v^{2,1}_n={(-1)^n e^{-nt_{2,1}} \over[n]}, \qquad 
\hat u^2_n={1\over [n]} -e^{-nt_{2,1}}u^2_n ={1 \over [n]} 
\biggl(1-\sum_{k=1}^{N_1} e^{-n(t_{2,1}+\cdots+ t_{2,k}) }\biggr).
$$
We conclude that
\begin{prop}
\begin{multline*}
Z_{N_1,N_2}(\la;\vt_1, \vt_2) 
  = \exp \biggl(F^1_{N_2}(\la;\vt_2)
  +F^2_{N_1}(\la;\vt_1)+ F^2_{N_2}(\la;\vt_2)\biggr) \\
  \sum_{\mu} q^{-\kappa_{\mu}/2} (-1)^{|\mu|} e^{-t_{1,1} |\mu|}
  s_{\mu}(\hat u^2)s_{\mu}(u^1),
\end{multline*}
where
\begin{eqnarray*}
 F^1_{N_i}(\la;\vt_i)&=& \sum_{n>0} {-1\over n[n]^2 }
            \sum_{k=1}^{N_i} e^{-n(t_{i,1}+\cdots + t_{i,k})}\\
 F^2_{N_i}(\la;\vt_i)&=& \sum_{n>0}\frac{1}{n[n]^2}
\sum_{2\leq k_1 \leq k_2\leq N_i} e^{-n(t_{i,k_1}+\cdots+t_{i,k_2})} 
\end{eqnarray*}
\end{prop}

\bibliographystyle{gtart}
\bibliography{link}

\end{document}